\documentclass{gtart_h}

\def\ifplaintex{\expandafter\ifx\csname documentclass\endcsname\relax}

\def\gtp{{\mathsurround=0pt\it $\cal G\mskip-2mu$eometry \&\ 
$\cal T\!\!$opology $\cal P\!$ublications}}  

\def\recd{{\small Received:\qua\receiveddate\ifx\reviseddate\relax
\else\qquad Revised:\qua\reviseddate\fi\par}} 

\def\lognumber#1{\def\thelognumber{#1}}
\def\volumenumber#1{\def\thevolumenumber{#1}}
\def\volumeyear#1{\def\thevolumeyear{#1}}
\def\papernumber#1{\def\thepapernumber{#1}}
\def\pagenumbers#1#2{\def\startpage{#1}\def\finishpage{#2}}
\def\published#1{\def\publishdate{#1}}

\def\received#1{\def\receiveddate{#1}}
\def\revised#1{\def\reviseddate{#1}}
\def\accepted#1{\def\accepteddate{#1}}
\def\asciititle#1{\def\theasciititle{#1}}
\def\covertitle#1{\def\thecovertitle{#1}}

\def\asciiaddress#1{\def\theasciiaddress{#1}}
\def\asciiemail#1{\def\theasciiemail{#1}}

\long\def\asciiabstract#1{\long\def\theasciiabstract{#1}}
\def\asciikeywords#1{\def\theasciikeywords{#1}}

\let\\\par\let\thelognumber\relax\let\thevolumenumber\relax
\let\thepapernumber\relax\let\thevolumeyear\relax\let\startpage\relax
\let\finishpage\relax\let\publishdate\relax\let\receiveddate\relax
\let\reviseddate\relax\let\accepteddate\relax\let\theasciititle\relax
\let\thecovertitle\relax\let\theasciiauthors\relax\let\theasciiaddress\relax
\let\theasciiabstract\relax\let\theasciikeywords\relax

\let\theasciiemail\relax

\ifplaintex
\font\logobig=cmssbx10 scaled 3836
\font\logomed=cmssbx10 scaled 2557
\else
\font\logobig=cmssbx10 scaled 4200
\font\logomed=cmssbx10 scaled 2800
\fi

\long\def\makeagttitle{   
\count0=\startpage
\agt\hfill      
\hbox to 45truept{\vbox to 0pt{\vglue -13truept{\logomed A\kern -.37em{\logobig 
T}\kern -.38em G}\vss}\hss}
\break
{\small Volume \thevolumenumber\ (\thevolumeyear)
\startpage--\finishpage\nl
Published: \publishdate}

\vglue .25truein

{\parskip=0pt\leftskip 0pt plus
1fil\def\\{\par\smallskip}{\Large\bf\thetitle}\par\medskip} \vglue
0.05truein

{\parskip=0pt\leftskip 0pt plus 1fil\def\\{\par}{\sc\theauthors}
\par\medskip}%
 
\vglue 0.03truein 

{\small\leftskip 25truept\rightskip 25truept{\bf Abstract}\stdspace\theabstract

{\bf AMS Classification}\stdspace\theprimaryclass
\ifx\thesecondaryclass\relax\else; \thesecondaryclass\fi\par
{\bf Keywords}\stdspace \thekeywords\par}\vglue 7truept

}   

\ifplaintex
\font\phead=cmsl9 scaled 950
\font\pnum=cmbx10 scaled 913
\font\pfoot=cmsl9 scaled 950
\headline{\vbox to 0pt{\vskip -4.5mm\line{\small\phead\ifnum
\count0=\startpage ISSN 1472-2739 (on-line) 1472-2747 (printed)
\hfill {\pnum\folio}\else\ifodd\count0\def\\{ }%
\ifx\theshorttitle\relax\thetitle\else\theshorttitle\fi\hfill{\pnum\folio}
\else\def\\{ and }{\pnum\folio}\hfill\ifx\theshortauthors\relax\theauthors
\else\theshortauthors\fi\fi\fi}\vss}}
\footline{\vbox to 0pt{\vglue 0mm\line{\small\pfoot\ifnum\count0=\startpage
\copyright\ \gtp\hfill\else
\agt, Volume \thevolumenumber\ (\thevolumeyear)\hfill\fi}\vss}}
\else
\font=cmsl9 scaled 1050
\font\lnum=cmbx10 
\font=cmsl9 scaled 1050
\makeatletter
\def\@oddhead{{\small\lhead\ifnum\count0=\startpage ISSN 1472-2739 
(on-line) 1472-2747 (printed)\hfill {\lnum\number\count0}\else\ifodd\count0
\def\\{ }\ifx\theshorttitle\relax \thetitle \else\theshorttitle\fi\hfill
{\lnum\number\count0}\else\def\\{ and }{\lnum\number\count0}
\hfill\ifx\theshortauthors\relax 
\theauthors\else\theshortauthors\fi\fi\fi}}\def\@evenhead{\@oddhead}
\def\@oddfoot{\small\lfoot\ifnum\count0=\startpage\copyright\ \gtp\hfill\else
\agt, Volume \thevolumenumber\ (\thevolumeyear)\hfill\fi}
\def\@evenfoot{\@oddfoot}
\makeatother
\fi
\let\maketitlepage\makeagttitle
\let\makeshorttitle\maketitlepage
\let\maketitle\maketitlepage

\newwrite\gtoutfile
\long\gdef\makeheadfile{  
{\def\\{, }\def\s{ }
\immediate\openout\gtoutfile head.xxx
\immediate\write\gtoutfile{Proxy-for: \ifx\theasciiauthors\relax
\theauthors\else\theasciiauthors\fi\s<\ifx\theasciiemail\relax\theemail\else\theasciiemail\fi>}
\immediate\write\gtoutfile{\noexpand\\}
\immediate\write\gtoutfile{Authors: \ifx\theasciiauthors\relax
\theauthors\else\theasciiauthors\fi}
{\def\\{ }\immediate\write\gtoutfile{Title: \ifx\theasciititle\relax
\thetitle\else\theasciititle\fi}}
\immediate\write\gtoutfile{Subj-class: GT or SG, GR etc}
\immediate\write\gtoutfile{MSC-class: \theprimaryclass\ifx\thesecondaryclass\relax\else, \thesecondaryclass\fi}
\immediate\write\gtoutfile{Journal-ref: Algebr. Geom. Topol. \thevolumenumber\s
(\thevolumeyear) \startpage-\finishpage}
\immediate\write\gtoutfile{Comments: Published by Algebraic and
Geometric Topology at}
\immediate\write\gtoutfile{\s\s\s  http://www.maths.warwick.ac.uk/agt/AGTVol\thevolumenumber/agt-\thevolumenumber-\thepapernumber.abs.html}
\immediate\write\gtoutfile{\noexpand\\}
\immediate\write\gtoutfile{}
\ifx\theasciiabstract\relax
\immediate\write\gtoutfile{\theabstract}\else
\immediate\write\gtoutfile{\theasciiabstract}\fi
\immediate\write\gtoutfile{}
\immediate\write\gtoutfile{\noexpand\\}
\immediate\write\gtoutfile{}
\immediate\closeout\gtoutfile}}  

\def\maketitlepage{\makeagttitle\makeheadfile}
\let\makeshorttitle\maketitlepage
\let\maketitle\maketitlepage

\lognumber{9}
\volumenumber{5}
\volumeyear{2005}
\papernumber{9}
\pagenumbers{135}{182}
\received{8 May 2003} 
\revised{18 August 2004}
\accepted{21 September 2004}
\published{12 March 2005}

\usepackage{amssymb,amsmath}
\usepackage[all]{xy}

\newcommand{\refequ}[1]{$(\ref{#1})$}

\newtheorem{thm}{Theorem}[section]
\newtheorem{corol}[thm]{Corollary}
\newtheorem{lemma}[thm]{Lemma}
\newtheorem{prop}[thm]{Proposition}

\theoremstyle{definition}
\newtheorem{defin}[thm]{Definition}
\newtheorem{example}[thm]{Example}
\newtheorem{examples}[thm]{Examples}
\theoremstyle{remark}
\newtheorem{rmk}[thm]{Remark}
\numberwithin{equation}{section}

\newcommand{\BR}{\mathbb R}
\newcommand{\BQ}{\mathbb Q}
\newcommand{\BD}{\mathbb D}
\newcommand{\BDt}{\underline{\mathbb D}}
\newcommand{\BZ}{\mathbb Z}
\newcommand{\Bk}{\mathbf {k}}
\newcommand{\calC}{\mathcal{C}}
\newcommand{\calS}{\mathcal{S}}

\newcommand{\phishrk}{{\psi}}
\newcommand{\phishrkk}{{\psi_k}}

\newcommand{\quism}{\stackrel{\simeq}{\rightarrow}}
\newcommand{\iso}{\stackrel{\cong}{\rightarrow}}

\newcommand{\topdegree}{top-degree\ }
\newcommand{\Topdegree}{Top-degree\ }
\newcommand{\Apl}{A_{PL}}
\newcommand{\homset}{{\mathrm{hom}}}
\newcommand{\id}{{\mathrm{id}}}
\newcommand{\isoset}{{\mathrm{iso}}}
\newcommand{\pr}{\mathrm{pr}}
\newcommand{\del}{\partial}

\newcommand{\set}[1]{\left\{{#1}\right\}}
\newcommand{\punct}{\check}
\newcommand{\takeaway}{\smallsetminus}

\begin{document}
\sloppy

\title{Algebraic models of Poincar\'e embeddings}
\asciititle{Algebraic models of Poincare embeddings}
\covertitle{Algebraic models of Poincar\noexpand\'e embeddings}
\author{Pascal Lambrechts\\Don Stanley}

\address{Institut Math\'ematique, Universit\'e de Louvain,
2, chemin du Cyclotron\\B-1348 Louvain-la-Neuve, Belgium}
\address{Department of Mathematics and Statistics, University of Regina\\
College West 307.14, Regina, Saskatchewan, Canada, S4S 0A2}

\asciiaddress{Institut Mathematique, Universite de Louvain,
2, chemin du Cyclotron\\B-1348 Louvain-la-Neuve, 
Belgium\\and\\Department of Mathematics and Statistics, University of Regina\\
College West 307.14, Regina, Saskatchewan, Canada, S4S 0A2}

\asciiemail{lambrechts@math.ucl.ac.be, stanley@math.uregina.ca}
\gtemail{\mailto{lambrechts@math.ucl.ac.be}{\rm\qua 
and\qua}\mailto{stanley@math.uregina.ca}}%


\asciikeywords{Poincare embeddings, Lefschetz duality, Sullivan models}

\begin{abstract}
Let $f\co P\hookrightarrow W$ be an embedding of a
compact polyhedron in a closed oriented
manifold $W$, let $T$ be a regular neighborhood of $P$ in $W$
and let $C:=\overline{W\takeaway T}$ be its complement.
Then $W$  is the homotopy push-out of 
a diagram $C\leftarrow \del T\to P$. This homotopy push-out square is
an example of what is called a Poincar\'e embedding. 

We  study how to construct  algebraic models, in
particular in the sense of Sullivan, of that homotopy push-out from
a model of the map $f$. When the codimension is high enough
this allows us to completely determine the rational homotopy type of the complement
$C\simeq W\takeaway f(P)$. Moreover we
construct examples to show that our restriction on the codimension
is sharp.

Without restriction on the codimension
we also give  differentiable modules models of Poincar\'e embeddings
 and we deduce a refinement
of the classical Lefschetz duality theorem, 
giving information on the algebra structure of the cohomology of the complement.
\end{abstract}
\asciiabstract{%
Let f: P-->W be an embedding of a compact polyhedron in a closed
oriented manifold W, let T be a regular neighborhood of P in W and let
C:=closure(W-T) be its complement.  Then W is the homotopy push-out of
a diagram C<--dT-->P. This homotopy push-out square is an example of
what is called a Poincare embedding.  We study how to construct
algebraic models, in particular in the sense of Sullivan, of that
homotopy push-out from a model of the map f. When the codimension is
high enough this allows us to completely determine the rational
homotopy type of the complement C = W-f(P).  Moreover we construct
examples to show that our restriction on the codimension is sharp.
Without restriction on the codimension we also give differentiable
modules models of Poincare embeddings and we deduce a refinement of
the classical Lefschetz duality theorem, giving information on the
algebra structure of the cohomology of the complement.}

\maketitle

\section{Introduction}
Let us recall the notion of a \emph{Poincar\'e embedding}:

\begin{defin}[\rm({Levitt \cite{Levitt}, and \cite[Section 5]{Klein} for a modern exposition})]
\label{def-Pemb}
Let $W$ be a Poincar\'e duality space of dimension $n$ and let
$P$ be a finite CW-complex of dimension $m$.
A \emph{Poincar\'e embedding} of $P$ in $W$ (of \emph{dimension} $n$ and \emph{codimension} $n-m$)
is a commutative diagram of
topological spaces
\begin{equation}
\label{diag-mainsquare} \xymatrix{ \del
T\ar@{->}[r]^i\ar@{->}[d]_k&
P\ar@{->}[d]^f\\
C\ar@{->}[r]_l&W }
\end{equation}
such that \refequ{diag-mainsquare} is a homotopy push-out,
$(P,\del T)$ and $(C,\del T)$ are Poincar\'e duality 
pairs\footnote{By abuse of terminology, by the \emph{pair} $(P,\partial T)$
we actually mean the pair $(P',\partial T)$ where $P'$ is the mapping cylinder of $i$,
and similarly for the pair $(C,\partial T)$}
 in dimension $n$,
and the map $i$ is $(n-m-1)$-connected.
\end{defin}

The motivating example of a Poincar\'e embedding arises when $W$ is a closed orientable PL-manifold of dimension $n$ and
 $f\co P\hookrightarrow W$ is a piecewise
linear  embedding of a compact polyhedron $P$ in $W$.
Alternatively we can also take $f$ to be a smooth embedding between smooth compact manifolds.
Then $f(P)$ admits a  regular neighborhood, that is 
 a codimension $0$ compact
submanifold $T\subset W$ that deformation retracts to $P$ (see
\cite[page 33]{RourkeSanderson}.) Let $C:=\overline{W\takeaway T}$
be the closure of the complement of $T$ in $W$. Then $C$ and $T$
are both compact manifolds of dimension $n$ with a common boundary
$\del T=\del C$ and $W=T\cup_{\del T} C$. The composition of the inclusion $\del T\hookrightarrow T$
with the retraction $T\quism P$ gives a map $i\co \del T\to P$ and we obtain the pushout \refequ{diag-mainsquare}.
If the polyhedron $P$ is of dimension $m$,
then a general position argument implies that the map $i$ is $(n-m-1)$-connected.
Of course $C$ has the homotopy type of the complement $W\takeaway f(P)$.

Thus morally a Poincar\'e  embedding is the homotopy generalization of a PL embedding.
Notice that, in Definition \ref{def-Pemb}, $\del T$ is just a topological space and not necessarily
a genuine boundary of a manifold $T$, and $W$ does not need to be a manifold.
Notice also  that
by a Poincar\'e embedding we mean \emph{all} of the diagram \refequ{diag-mainsquare} and not only
the map $f$. When such a diagram exists we say that the map $f\co P\to W$ \emph{Poincar\'e embeds}.
The space $C$ in the push-out diagram is called the \emph{complement} of $P$.

A natural question is whether the {homotopy class} of a map $f$ that Poincar\'e embeds 
determines the square \refequ{diag-mainsquare} up to homotopy
equivalence and in particular the homotopy type of the complement $C$.
The answer to this question is negative in general as it can be seen
with $W=S^3$ and $P=S^1$. Indeed all PL-embeddings $f\co S^1\hookrightarrow S^3$ 
are nullhomotopic but the homotopy type of the complement $C\simeq S^3\takeaway f(S^1)$
can vary considerably (see for example
\cite[Corollary 11.3]{Lickorish} or \cite{GordonLuecke}.)
This is possible since in general the homotopy class $[f]$ of $f$ does not determine its
isotopy class. On the other hand in the case of a PL-embedding when the codimension is high enough, namely when $n\geq 2m+3$, 
then a general position argument  implies that $[f]$ determines the isotopy class of $f$. Therefore under
this high codimension hypothesis the homotopy class of a PL-embedding $f$ \emph{does} determine 
the homotopy type of the square \refequ{diag-mainsquare}. Similarly under a slightly more restrictive condition
on the codimension, there exists a unique Poincar\'e embedding \refequ{diag-mainsquare} associated to a given homotopy
class $[f]$. See Theorem \ref{thm-unknot} below
for a precise and more general statement for PL-embeddings as 
well as a discussion on the corresponding result for Poincar\'e embeddings.

The aim of this paper is to study an algebraic translation of the above question:
can we build  algebraic models, such as Sullivan models which
encode rational homotopy type, of the square \refequ{diag-mainsquare} from an algebraic model of
the map $f\,$? In order to be more precise, we first review Sullivan's theory for modeling
rational homotopy types by algebraic models. By a \emph{CDGA},
$A$, we mean a non-negatively graded algebra over the field $\BQ$ of
rational numbers that is commutative in the graded sense and
endowed with a degree $+1$ derivation $d\co A\to A$ such that
$d^2=0$. Sullivan has defined in \cite{Sullivan} a contravariant
functor from topological spaces to CDGA,
$$\Apl\co\mathrm{Top}\to\mathrm{CDGA},$$
mimicking the de Rham complex of differential forms on a smooth
manifold.
By a \emph{CDGA model} of a space $X$, we mean a CDGA,
$A$, linked to $\Apl(X)$ by a chain of CDGA morphisms inducing
isomorphisms in cohomology,
$$
\xymatrix{A&\ar[l]_\simeq A_1\ar[r]^\simeq&\cdots&\ar[l]_\simeq
A_n\ar[r]^-\simeq&\Apl(X).}
$$
The fundamental result of Sullivan's theory is that if $X$ is a
simply-connected space with rational homology of finite type,
then any  CDGA model of $X$ determines its rational homotopy type.
There is a similar result for maps and
 more generally for finite diagrams. See \cite{FHT-RHT} for a complete
exposition of that theory.

Our first result is the construction, under the high codimension hypothesis
 $\dim(W)\geq 2\dim(P)+3$, of an explicit CDGA model
of the Poincar\'e embedding \refequ{diag-mainsquare} out of a CDGA-model of $f$.
 To explain this result, we  need
some notation which will be made more precise in  Section 2.
We denote by $\# V:=\hom(V,\Bk)$ the dual of a $\Bk$-vector space $V$ and by $s^{p}X$
the $p$-th suspension of a graded object $X$, i.e.\ $(s^{p}X)^k=X^{p+k}$. The mapping cone of a cochain map
$f\co M\to N$ is written $N\oplus_fsM$. When $N$ is a CDGA and $M$ is an $N$-DGmodule
this mapping cone
can be endowed with  the multiplication 
$(n,sm)\cdot(n',sm')=(n\cdot n',s(n\cdot m'\pm n'\cdot m))$.
The differential of the mapping cone does not always satisfy the Leibnitz rule
for this multiplication, but it does under certain conditions on the dimensions 
and then the induced structure is called the
\emph{semi-trivial CDGA-structure} on the mapping cone (Definition \ref{def-CDGAMC}).

Our goal is to build 
 a CDGA model of the homotopy push-out \refequ{diag-mainsquare}, and in particular of the complement $C$,
 out of a CDGA model $\phi\co R\to Q$ of $f^*\co \Apl(W)\to\Apl(P)$.
Motivated by Lefschetz duality a first guess for a model of $\Apl(C)$ is the mapping cone
 $$
 R\oplus_\psi ss^{-n}\#Q
 $$
 $$\psi\co s^{-n}\#Q\to R\leqno{\rm where}$$
 is an $R$-DGmodule map such that $H^n(\psi)$ is an isomorphism.
 Unfortunately this naive guess has two flaws:
 \begin{enumerate}
 \item[(A)] such a map $\psi$ does not necessarily exist, and
 \item[(B)] the multiplication on $R\oplus_\psi ss^{-n}\#Q$ does not necessarily define a CDGA
 structure  because of the possible failure of the Leibnitz rule.
 \end{enumerate}
Problem (A) can be addressed by replacing $s^{-n}\#Q$ by a suitable weakly equivalent DG-module $D$, for example a cofibrant one, for
which there exists a map $\psi\co D\to R$ inducing an isomorphism in cohomology in degree $n$. Such a map is called 
a \emph{\topdegree map}\footnote{It was called
a \emph{shriek map} in earlier versions of this paper.} in Definition \ref{def-shriek}.
Problem (B) can be solved by restricting the range of degrees of the graded objects $R$, $Q$, and $D$. This is where
the high codimension hypothesis is needed.

We can now state our first result:
\begin{thm}\label{thm-stableCDGA}%
Consider a Poincar\'e embedding \refequ{diag-mainsquare}
with $P$ and $W$ connected. 
If $n\geq2m+3$ and $H^1(f;\BQ)$ is injective
then a model of the commutative CDGA square
$$\BD':=\vcenter{
\xymatrix@1{%
\Apl(W)\ar[r]^{f^*}\ar[d]_{l^*}&
\Apl(P)\ar[d]^{i^*}\\
\Apl(C)\ar[r]_{k^*}&%
\Apl(\del T) }}
$$
can be build explicitly out of any CDGA model of $f^*\co\Apl(W)\to\Apl(P)$.

More precisely, if $n\geq 2m+4$ or if $n\geq2m+3$ and $H^1(f;\BQ)$ is injective, then
the commutative CDGA square $\BD'$ is weakly equivalent to any commutative CDGA square
$$
\BD:=\vcenter{
\xymatrix{%
 R\ar[r]^\phi\ar@{^(->}[d]&%
Q\ar@{^(->}[d]\\
R\oplus_{\phishrk}sD\ar[r]_{\phi\oplus\id}&%
 Q\oplus_{}sD%
}}$$
where 
\begin{enumerate}
\item[\rm(i)] $\phi\co R\to Q$
is a CDGA model of $f^*:\Apl(W)\to\Apl(P)$ with $R^{>n}=0$ and $Q^{>m+2}=0$;
\item[\rm(ii)] $D$ is a $Q$-DGmodule weakly equivalent to $s^{-n}\#Q$ with $D^{>n+1}=0$ and $D^{<n-m}=0$;
\item[\rm(iii)] $\phishrk\co D\to R$ is an $R$-DGmodules map such that $H^n(\phishrk)$ is an isomorphism
\end{enumerate}
and the mapping cones are endowed with the semi-trivial CDGA structure.

Moreover if $n\geq 2m+3$ and
$H^1(f;\BQ)$ is injective, then $R$, $Q$, $D$, $\phi$, and $\psi$ 
satisfying (i)-(iii) can be \emph{explicitly} constructed  out of any CDGA model of
$f^*\co \Apl(W)\to\Apl(P)$.
\end{thm}

Since CDGA models encode rational homotopy types of simply connected spaces
an immediate corollary of the above theorem is that when $P$ and $W$ are simply connected and 
$\dim(W)\geq 2\dim(P)+3$, then the rational homotopy type of the Poincar\'e  embedding
\refequ{diag-mainsquare} depends only on the rational homotopy class of $f$.

As a byproduct of this theorem  we obtain also a CDGA model $Q\oplus ss^{-n}\#Q$ of
the boundary $\del T$ of a thickening of $P$ under a high codimension hypothesis. This model was already
described in \cite{Lambrechts-thickening} and an analogous model is built in \cite{KahlLVdb}
under weaker hypotheses.

\medbreak

In our first theorem we have supposed that $\dim W\geq2\dim P+3$.
When the connectivity of the embedding is high this condition on the
codimension can be weakened. Indeed in the case of PL-embeddings we
 have the following classical result:

\begin{thm}[PL-unknotting, Wall and Hudson]\label{thm-unknot}
Let $P$ be a compact $m$-dimensional polyhedron and let $W$
be a closed $n$-dimensional manifold with $n\geq m+3$. Let $r$ be an integer such that
\begin{equation}\label{equ-unknot}
r\geq2m-n+2.
\end{equation}
Then any two homotopic $r$-connected embeddings  $f_0,f_1\co
P\hookrightarrow W$ are isotopic. As a consequence, if $f$ is $r$-connected then the homotopy type
of the square \refequ{diag-mainsquare} depends only on the homotopy class of $f$.
\end{thm}
\begin{proof}
By the uniqueness
part of the Wall's embedding theorem \cite[page 76]{Wall-thick}
$f_0$ and $f_1$ are concordant. Since the codimension
is at least $3$, concordance implies isotopy
\cite{Hudson-conc=>iso}. Therefore $f_0$ is isotopic to $f_1$.
By the uniqueness of a regular neighborhood this implies that the squares \refequ{diag-mainsquare}
for $f_0$ and $f_1$ are homeomorphic.
\end{proof}
The hypothesis that $f$ is $r$-connected with $r$ satisfying the inequality
\refequ{equ-unknot} is called the \emph{unknotting condition}.
The reason for which we have stated Theorem \ref{thm-unknot} in the context of PL-embeddings instead of Poincar\'e 
embeddings is that the corresponding result
for Poincar\'e embeddings is known only under a slightly more restrictive condition.
Indeed Klein has proved such an uniqueness result for Poincar\'e embeddings
with an unknotting condition increased by one, i.e.\ $r\geq2m-n+3$ \cite[Theorem 5.4]{Klein}, or
with the sharp unknotting condition \refequ{equ-unknot} 
 in the metastable range \cite{Klein-compression}. It  is still an open question
  whether 
condition \refequ{equ-unknot} guarantees the uniqueness
of Poincar\'e embeddings in full generality. 

We will prove a rational homotopy
theoretical partial  version of Theorem \ref{thm-unknot} by establishing
that, under the unknotting condition \refequ{equ-unknot}, the
rational homotopy type of the complement $C$ depends
only on the rational homotopy class of $f$. From Theorem \ref{thm-stableCDGA}
a guess for the model of the complement would be $R\oplus_\psi sD$ with
some assumption on the vanishing of $R$, $Q$, and $D$ in high degrees.
This vanishing assumption can be removed if we truncate the mapping cone $R\oplus_\psi sD$ by a suitable
acyclic module $L$. Moreover only a structure of $R$-DGmodule (instead of $Q$-DGmodule) is needed on
$D$. More precisely we have the following theorem:
\begin{thm}\label{thm-wkstCDGA}
Consider a Poincar\'e embedding \refequ{diag-mainsquare} of codimension
at least $2$ with $P$ and $W$ connected.
Let $r$ be a positive integer such that 
$H_*(f;\BQ)$ is $r$-connected, that is $H_i(f;\BQ)$ is an isomorphism for $i<r$ and
an epimorphism for $i=r$.

If
\begin{equation}\label{equ-unknotrht}r\geq2m-n+2.\end{equation}
then a CDGA model of the map $l\co C\to W$ can be build explicitly out of
any CDGA model of $f\co P\to W$.

More precisely, let
\begin{enumerate}
\item[\rm(i)] $\phi\co R\to Q$ be a CDGA model of
$f^*\co\Apl(W)\to\Apl(P)$ with $R$  connected;
\item[\rm(ii)] $D$ be an $R$-DGmodule weakly equivalent to $s^{-n}\#Q$ with
$D^{<n-m}=0$;
\item[\rm(iii)] $\phishrk\co D\to R$ be a \topdegree map of $R$-DGmodules;
\item[\rm(iv)] $L\subset R\oplus_{\phishrk}sD$ be an acyclic $R$-subDGmodule 
with $L^{\leq n-r-2}=0$ and $(R\oplus_{\phishrk}sD)^{\geq n-r}\subset L$.
\end{enumerate}
Then the canonical CDGA map
$$\lambda\co R\to (R\oplus_{\phishrk}sD)/L$$
is a CDGA-model of the map
$$l^*\co\Apl(W)\to\Apl(C).$$
where  $\lambda$ is the composition of the inclusion with the projection and the algebra structure
on the truncated mapping cone is induced by the formula
$(r,sd)\cdot(r',sd')=(r\cdot r',s(r\cdot d'\pm r'\cdot d))$.

Moreover under condition \refequ{equ-unknotrht} it is possible to construct explicitly $R$, $Q$, $D$, $L$,
$\phi$, $\psi$ satisfying hypotheses (i)--(iv) out of any CDGA-model of $f^*\co\Apl(W)\to\Apl(P)$.
\end{thm}

\begin{corol}
\label{corol-wkstCDGA} Consider a Poincar\'e embedding \refequ{diag-mainsquare}
of codimension at least $3$ and with $P$ and $W$ simply-connected.
Let $r$ be a positive integer such that $H_*(f;\BQ)$ is $r$-connected. If $r\geq 2m-n+2$
then the rational homotopy type of the complement $C$ depends only on the rational homotopy class of $f$.
\end{corol}
Moreover we will show that the unknotting condition in Theorem \ref{thm-wkstCDGA} 
is sharp. More precisely we will construct in Propositions \ref{prop-exsharp1} and \ref{prop-exsharp2} 
families of examples for which the unknotting condition \refequ{equ-unknotrht} fails only by a little
but such that the rational cohomology algebra of the complement is not determined
by the rational homotopy class
of the embedding. Note also that our rational result is valid for any 
Poincar\'e embeddings satisfying the unknotting condition,
which improves by $1$ the hypothesis under which the ``integral''
homotopy type of the complement
is known to be unique \cite[Corollary B]{Klein-2}.

Unfortunately we were not able to determine the complete rational homotopy type of the square
\refequ{diag-mainsquare} from the rational homotopy class of $f$ under the unknotting condition.
The best result that we can prove in this direction is the determination, under connectivity
hypotheses
on $P$ and $W$ and the extra assumption that
$n\geq m+r+2$, of the modified
square \refequ{diag-mainsquare} where $\del T$ is replaced by the space
$\punct{\del T}$ obtained by removing its top cell. See Theorem \ref{thm-wkstCDGAsquare} for
a precise statement.

Our rational models in Theorems \ref{thm-stableCDGA} and \ref{thm-wkstCDGA}
have applications to the construction of the model of blow-ups \cite{LS-stableblowup} and
\cite{LS-unstableblowup}, and of
the configuration space on two points  \cite{LS-FM2}.
\medbreak

The above discussion was about CDGA models for the square \refequ{diag-mainsquare}
 which determine its rational homotopy type. Instead of CDGA models  associated to the functor
$\Apl$
we can associate models to the functor of singular cochains with
coefficients in a field $\Bk$ of arbitrary characteristic,
$S^*(-;\Bk)$. If $Y$ is a space then $S^*(Y;\Bk)$ is a
differential graded algebra (a DGA for short), and if $f\co
X\to Y$ is a continuous map then $S^*(X;\Bk)$ is a differential
graded module (\emph{DGmodule}) over the DGA $S^*(Y;\Bk)$. There
is a notion of models of such DGmodules, and we can build such a model of
 the Poincar\'e embedding \refequ{diag-mainsquare} without any restriction on the
codimension or even on the connectivity of $P$. To state the result we use the notion of a
\emph{menorah} as defined in Example \ref{examples-diagram} and
which is essentially a family of maps with same domain.
\begin{thm}\label{thm-DGmodnonconn}%
Consider a Poincar\'e embedding  \refequ{diag-mainsquare} with $W$ connected.
Denote the connected components of $P$ by $P_1,\cdots,P_c$  and
set $f_k:=f|P_k$, for $k=1,\cdots,c$. Denote by $C^*$ one of the
functors $S^*(-;\Bk)$ or $\Apl$.

Suppose a quasi-isomorphism of DGA $\rho\co A\quism
C^*(W)$ has been given. Let
$$\set{\phi_k\co R\to Q_k}_{1\leq k\leq c}$$
be a model in $A$-DGMod of the menorah
$$\set{C^*(f_k)\co C^*(W)\to C^*(P_k)}_{1\leq k\leq c}.$$
For $k=1,\cdots,c$, let $D_k$ be an $A$-DGmodule weakly equivalent
to $s^{-n}\#C^*(P_k)$ and let
$
\phishrk_k\co D_k\to R$ be a \topdegree map of $A$-DGmodules.

Set $D=\oplus_{k=1}^cD_k$,  $Q=\oplus_{k=1}^cQ_k$,
$\phi=(\phi_1,\ldots,\phi_c)\co R\to Q$, and
$\phishrk=\sum_{k=1}^c\phishrk_k\co D\to R$. Then the two following
commutative
squares are weakly equivalent in $A$-DGMod:
$$\BD:=\vcenter{
\xymatrix@1{%
 R\vrule width0pt depth6pt\ar[r]^\phi\ar@{^(->}[d]&%
Q\vrule width0pt depth6pt\ar@{^(->}[d]\\
R\oplus_{\phishrk}sD\ar[r]_{\phi\oplus\id}&%
 Q\oplus_{\phi\phishrk}sD%
}}\hbox{\quad\quad and \quad\quad}\BD':=\vcenter{
\xymatrix@1{%
C^*(W)\ar[r]^{f^*}\ar[d]_{l^*}&
C^*(P)\ar[d]^{i^*}\\
C^*(C)\ar[r]_{k^*}&%
 C^*(\del T). }}
$$
\end{thm}

This DGmodule model enables us  to improve the classical Lefschetz duality
theorem. Indeed this classical result states that the cohomology of the complement,
$H^*(C;\Bk)=H^*(W\takeaway f(P);\Bk)$, is determined \emph{as a vector space}
by the algebra map $H^*(f)\co  H^*(W)\to H^*(P)$.
Our result gives a way to determine the \emph{$H^*(W)$-module structure} of
$H^*(C)$, and even its algebra structure under the unknotting condition. This is the content 
of the following:

\begin{corol}[Improved Lefschetz duality]
\label{corol-HWmodstruct}
Consider  a Poincar\'e embedding \refequ{diag-mainsquare} with  $W$
connected.
Suppose a quasi-isomorphism of DGA \,$\rho\co A\quism
C^*(W)$ has been given and let $\phi\co R\to Q$ be an  $A$-DGmodule model of
$f^*\co C^*(W)\to C^*(P)$. Then we have an
isomorphism of $H^*(W;\Bk)$-modules
$$H^*(C;\Bk)\cong H(s^{-n}\#R\oplus_{s^{-n}\#\phi}s(s^{-n}\#Q)).$$
If moreover
$H_*(f;\Bk)$ is $r$-connected with $r\geq 2m-n+2$
then this isomorphism determines the algebra structure on $H^*(C;\Bk)$.
\end{corol}
Examples of Section 9 will show that the unknotting condition cannot be dropped
when determining the algebra structure in the last corollary. 
Christophe Boilley \cite{Boilley}
has constructed examples showing that the $H^*(W)$-module
structure on $H^*(C)$ is neither necessarily given by a trivial extension
nor determined by the map $H^*(f)$ induced in cohomology.

\medbreak

Notice that in all the results of this paper we can replace the Poincar\'e
embedding by the following weaker notion. Let  $\Bk$
be a field. A \emph{$\Bk$-Poincar\'e embedding}
 is a commutative square \refequ{diag-mainsquare} such that $W$, $(P,\del T)$ and $(C,\del T)$
satisfy Poincar\'e duality in dimension $n$ over $\Bk$,
$m$ is the cohomological dimension of $P$ with coefficients in $\Bk$,
$H^*(i;\Bk)$ is $(n-m-1)$-connected,
and the square \refequ{diag-mainsquare}  induces a Mayer-Vietoris long exact sequence in $H^*(-;\Bk)$.
In other words such a $\Bk$-Poincar\'e embedding is a \emph{homological} version of a
Poincar\'e embedding.

As a last remark note that our study is 
complementary to the work of Morgan
\cite{Morgan} 
who has computed the rational homotopy type of the complement of
divisors $D_i$ with normal crossings in a projective algebraic
variety $W$. In his case the codimension is very low ($D_i$ is of
codimension $2$) but the existence of mixed Hodge structures
 \cite{Deligne}
implies that the rational homotopy type of the complement is
determined by the maps induced in cohomology by the inclusion of
divisors. In the case of a single divisor $D$,
 Morgan's model for $W\takeaway D$ is
expressed in terms of the  shriek map
$f^!\co H^{*+2}(D)\to H^*(W)$ which is a special case of our \topdegree map
(see Example \ref{ex-topdegreeshriek}.)

\medskip

{\bf Plan of the rest of the paper}\qua Section 2 contains notation
and terminology and Section 3 is about diagrams in closed model
categories.  We explain in this section what we mean by a model of a
square or a menorah.  In Section 4 we define the notion of a
semi-trivial CDGA structure on certain mapping cones and in Section 5
we study the notion of a \topdegree map and prove their existence and
essential uniqueness.  Section 6 is about the DGmodule model of a
Poincar\'e embedding and contains the proofs of Theorem
\ref{thm-DGmodnonconn} and Corollary \ref{corol-HWmodstruct}.  Section
7 is about CDGA models of a Poincar\'e embedding in the stable case
and contains the proof of Theorem \ref{thm-stableCDGA}.  Section 8
discusses CDGA models of the complement in a Poincar\'e embedding
under the unknotting condition.  We prove here Theorem
\ref{thm-wkstCDGA} and its corollaries.  We also state and prove
Theorem \ref{thm-wkstCDGAsquare} which exhibits a model of a square
related to \refequ{diag-mainsquare} under a stronger unknotting
condition.  Finally Section 9 contains examples of rationally knotted
embeddings and we illustrate by explicit examples the sharpness of the
unknotting condition.

\medskip{\bf Acknowledgements}\qua The authors want to thank Bill Dwyer
for enlightening conversations on closed model structures on
categories of diagrams and John Klein for explaining the proof of
Theorem \ref{thm-unknot}. We thank also the referee for pointing out
that our results could apply to Poincar\'e embeddings.  During this
work the first author benefited from the hospitality of the University
of Alberta and of a travel grant from F.N.R.S., and the second author
from the hospitality of the Universit\'e of Louvain. The first author
is Chercheur Qualifi\'e au F.N.R.S.

\section{Notation and terminology} \label{section-toolkitApl}We
denote by $\Bk$ a commutative field. Recall the notions of
\emph{differential graded algebra}, or DGA for short, and of
(left) \emph{graded differential modules} over a DGA $R$, or
$R$-DGmodules for short, as both defined for example in
\cite[Section 3(c)]{FHT-RHT}. We will always suppose that the
DGA are non negatively graded and that the
differentials are of degree $+1$. We denote by $R$-DGMod the
category of $R$-DGmodules.

{\bf Convention on left and right modules}\qua Sometimes in the paper
(in particular in Section \ref{section-DGMod}) it will be important to
distinguish between left and right DGmodules. By an
\emph{$R$-DGmodule} we always mean a {left} $R$-DGmodule, otherwise we
write explicitly \emph{right $R$-DGmodule}. Also by $R$-DGMod we
denote only the category of left $R$-DGmodules. We denote by
$\hom_\Bk$ (resp.\ $\hom_R$) the sets of $\Bk$-modules (resp.\
$R$-modules) morphisms.
 
 We have also a notion of
\emph{commutative differential graded algebra}, or CDGA for short,
which is a DGA such that the multiplication is graded commutative
(\cite[Example 5 in Section 3(b)]{FHT-RHT} where there are  called
\emph{commutative cochain algebras}). We denote
by CDGA the corresponding category.  A CDGA or more generally a
non-negatively graded vector space, $V$, is called
\emph{connected} if $V^0\cong\Bk$.

The degrees of graded modules and algebras will be written as
superscripts. If $X$ is a graded module or algebra, we will write
$X^{>m}=0$ to express the fact that $X^k=0$ for $k>m$, and
similarly $X^{\geq m}=0$, $X^{<m}=0$, and so on.

The \emph{dual} of a graded $\Bk$-module $M$ will be
denoted by $\#M$ with the grading $ (\#M)^i =
\homset(M^{-i},\Bk)$. The duality pairing is defined by
$$
\langle-,-\rangle\co M\otimes \#M \to\Bk,\,x\otimes f\mapsto
\langle x,f\rangle=f(x).
$$
If $(M,d)$ is a differential module then its dual $\#M$ is
equipped with the differential $\delta$ characterized by
$\langle x, \delta(f)\rangle=-(-1)^{|x|}\langle d(x),f\rangle$.
If $M$ is a \emph{right} module over some graded algebra $R$, then
its dual admits a structure of \emph{left} $R$-module
characterized by the formula $ \langle x,a.f \rangle=\langle x.a,f
\rangle$. Similarly if $M$ is a right DGmodule then its dual
becomes a left DGmodule.

The {\em $k$-th suspension} of a graded vector space $M$ is the
graded vector space $s^kM$ defined by $(s^kM)^j\cong M^{k+j}$ and
this isomorphism is denoted by $s^k$. If $M$ is also a left
$R$-module, we transport this structure on  the $k$-suspension by
the formula $r.(s^kx)=(-1)^{|r|k}s^k(r.x)$. Also if $M$ is
equipped with a differential $d$, then we define a differential on
$s^kM$ by $d(s^kx)=(-1)^ks^k(dx)$. If $k=1$ we write $sM$ for
$s^1M$.

The {\em mapping cone} of an $R$-DGmodule morphism $f\co X\to Y$
is the $R$-DGmodule $C(f):=(Y\oplus_f sX,d)$ where the
differential is defined by $d(y,sx)=(d_Y(y)+f(x),-sd_X(x))$. If
$f$ is a CDGA morphism, in general there is no natural CDGA
structure on the mapping cone but we will show in Section
\ref{section-MC} that such a CDGA structure exists under favorable
hypotheses.

We will use the functor of (normalized) singular cochains with
coefficients in $\Bk$
$
S^*(-;\Bk)\co\mathrm{ Top}\to \mathrm{DGA}
$
as defined for example in \cite[Chapter 5]{FHT-RHT}. When $\Bk$ is
of characteristic $0$, we have also the de Rham-Sullivan functor
of polynomial forms
$
\Apl\co \mathrm{Top}\to \mathrm{CDGA}
$
as defined in \cite{BG} or \cite[Chapter 10]{FHT-RHT}.

The categories $R$-DGMod and CDGA are closed model categories in the
sense of Quillen for which the weak equivalences are the
quasi-isomorphisms and the fibrations are the surjections (for a
nice review of closed model categories, we refer the reader to
\cite{DwyerSpalinski}). By an \emph{acyclic (co)fibration} we mean
a (co)fibration that is also a weak equivalence. We say that two
objects $X$ and $X'$ in a closed model category are \emph{weakly
equivalent} or that $X$ is a \emph{model} of $X'$ if there exists
a finite chain of weak equivalences joining them,
$$
\xymatrix{X&\ar[l]_\simeq X_1\ar[r]^\simeq&\cdots&\ar[l]_\simeq
X_n\ar[r]^\simeq&X'}.
$$
In that case we will write $X\simeq X'$. Since in Section
\ref{section-diagram} we will consider a closed model structure on
certain categories of diagrams, we can speak of \emph{models} of
that diagrams.

We review quickly the notion of \emph{relative Sullivan algebras}
which is an important class of cofibrations in CDGA. If $V$ is a
non-negatively graded vector space we denote by $\wedge V$ the
free graded commutative algebra generated by $V$ (see \cite[\S
3(b), Example 6]{FHT-RHT}.) A {relative Sullivan algebra}
(\cite[Chapter 14]{FHT-RHT}, or \emph{KS-extension} in the older
terminology of \cite{Halperin-lecturesminmod}) is a CDGA morphism
$\iota\co(A,d_A)\hookrightarrow (A\otimes \wedge V,D)$
where the differential $D$ is an extension of $d_A$ that satisfies
some nilpotence condition (see \cite[Chapter 14]{FHT-RHT} for the
precise definition.) Notice that in this paper we do not assume
that $V^0=0$, following  \cite{Halperin-lecturesminmod} but
contrary to \cite{FHT-RHT}. In the special case $A=\Bk$ we get the
notion of a \emph{Sullivan algebra}, $(\wedge V,D)$, which is a
cofibrant object in CDGA. Examples of cofibrant objects in
$R$-DGMod are \emph{semi-free models} as defined in \cite[Chapter
6]{FHT-RHT}. Roughly speaking they are $R$-DGmodules of the form
$(R\otimes V,D)$ where $V$ is a graded vector space and the
differential $D$ satisfies also a nilpotence condition. Finally
remember that every object is fibrant in CDGA and in $R$-DGMod.

To denote that two maps $f_0$ and $f_1$ are homotopic in CDGA or
$R$-DGMod  we will write $f_0\sim f_1$, or sometimes $f_0\sim_R
f_1$ to emphasize the underlying DGA. When $P$ and $N$ are
$R$-DGmodules, with $P$ cofibrant, we denote by
$$
[P,N]_R
$$
the set of homotopy classes of $R$-DGmodules from $P$ to $N$.

\section{Diagrams in closed model categories}\label{section-diagram}
In order of being able to speak of models
of objects, maps, commutative squares, and so on, we review in
this section the convenient language of diagrams as described for example in
\cite[Section 10]{DwyerSpalinski}. There will  exist a
closed model structure on each of the categories of diagrams that we will
consider. We will finish the section by two useful lemmas to turn
certain homotopy commutative diagrams into commutative ones.
\begin{defin}
Let $\calS$ be a small category and let $\calC$ be any category. A
\emph{diagram in $\calC$ shaped on $\calS$} is a covariant functor
$\BD\co\calS\to\calC$ and we say that $\calS$ is \emph{shaping} the diagram.
A \emph{morphism of diagrams} is a natural transformation between
two diagrams. This defines the category of diagrams $\calC^\calS$.
\end{defin}
We describe now the five main examples of diagrams that we will
consider in this paper. First recall that to each partially
ordered set (or \emph{poset}, for short), $({S},\leq)$, we can
associate a small category $\calS$ whose objects are the elements
of ${S}$ and such that the set of morphisms, $\homset_\calS(x,y)$,
between two objects $x$ and $y$ in $\calS$ is a singleton if
$x\leq y$ and is the empty set otherwise.
\begin{examples}$\phantom{999}$
\label{examples-diagram}
\begin{description} 
  \item[Object] If $\calS$ is the category with only one object and one morphism
  (that is the category associated with the poset with only one element) then a diagram
  in $\calC$ shaped on $\calS$ is called \emph{an object} of $\calC$.
  \item[Map]If $\calS$ is the category associated to the ordered set $\set{0,1}$
  then a diagram
  in $\calC$ shaped on $\calS$ is just a map between two objects of $\calC$.
  Such a diagram is called \emph{a map} of $\calC$.
  \item[Commutative square] Let $\calS$ be the category whose objects are  the four sets $\emptyset$, $\set{1}$,
   $\set{2}$, and $\set{1,2}$,
  and whose morphisms are the inclusion maps. A diagram in $\calC$ shaped on $\calS$ is called a
  \emph{commutative square} in $\calC$.
  \item[Menorah] Let $\calS$ be the category whose objects are
  $\emptyset,\set{1},\cdots,\set{n}$, for some positive integer $n$
  and where morphisms are inclusions of sets. 
  Then a diagram in $\calC$ shaped on $\calS$  is just a collection of maps $f_1,\cdots,f_n$ with same domain.
  We call such a diagram \emph{a menorah
}
  and we denote it by $\set{f_i}_{1\leq i\leq
  n}$.
  \item[Composite] Let $\calS$ be the category corresponding to the ordered
  set $\set{0,1,2}$. A diagram shaped on $\calS$ is just two
  composable maps $f_0\co X\to Y$ and $f_1\co Y\to Z$. We call such a diagram a \emph{composite}
  and  we denote it by
  $(f_0,f_1)$.
\end{description}
\end{examples}
Each category shaping one of the five diagrams in Example
\ref{examples-diagram} is a \emph{very small category} in the
sense of \cite[Section 10.13]{DwyerSpalinski}. This notion is
useful because of the following:
\begin{prop}\label{prop-CMdiagrams}
Let $\calC$ be a closed model category and let $\calS$ be a very
small category. Then the category $\calC^\calS$ of diagrams in
$\calC$ shaped on $\calS$ admits a closed model structure such
that a map $f\co \BD\to \BD'$ between diagrams is a weak
equivalence (resp.\ a fibration) if and only if for each object $x$
in $\calS$ the map $f(x)\co \BD(x)\to \BD'(x)$ is a weak
equivalence (resp.\ a fibration) in $\calC$.

Moreover if $\hat \BD$ is a cofibrant diagram in $\calC^\calS$ then
for each object $x$ in $\calS$, $\hat \BD(x)$ is a cofibrant
object of $\calC$, and for each morphism $i$ in $\calS$, the map
$\hat \BD(i)$ is a cofibration in $\calC$.

If every object of $\calC$ is fibrant, then the same is true in
$\calC^\calS$.
\end{prop}
\begin{proof} This
model structure is described  in \cite[Section
10.13]{DwyerSpalinski}, where the cofibrations in $\calC^\calS$
are also defined (a complete proof of the axioms of Quillen for
this category can be found in \cite[Theorem 5.2.5]{Hovey}). Using
the fact that the initial object $\emptyset$ in $\calC^\calS$ is
the constant diagram with value $\emptyset$ at each object of
$\calS$, it is straightforward to check from the definition of a
cofibration in $\calC^\calS$ (\cite[10.13]{DwyerSpalinski}) that
if $\emptyset\to \hat \BD$ is a cofibration  then each object
$\hat \BD(x)$ is cofibrant and each map $\hat \BD(i)$ is a
cofibration. The last statement is obvious.
\end{proof}

In this paper we will always suppose that the closed model
structure on a category of diagrams $\calC^\calS$ is the one
considered in Proposition \ref{prop-CMdiagrams}. 
Following the terminology of Section \ref{section-toolkitApl} we can speak of weakly
equivalent diagrams or of a model of a diagram.

\begin{rmk}
\label{rmk-modelmenorah} If a menorah $\{f_k\}_{1\leq k\leq n}$
is a model of another menorah $\{f'_k\}_{1\leq k\leq n}$, then
clearly each map $f_k$ is a model of $f'_k$. It is important to
notice that the converse is \emph{not} true in general. Similarly
if a composite $(f,g)$ is a model of a composite $(f',g')$ then
$f$ is a model of $f'$ and $g$ is a model of $g'$, but again the
converse is not true.
\end{rmk}

The proofs of the following two lemmas are based on standard
techniques of closed model categories and we leave them as
exercises for the reader.
\begin{lemma}\label{lemma-bisurj}%
Let $X$ and $X'$ be two weakly equivalent objects in some closed
model category in which every object is fibrant. Then there exists
a cofibrant object $\hat X$ and acyclic fibrations
$$
\xymatrix{X&\ar@{->>}[l]_\simeq^\beta \hat
X\ar@{->>}[r]^\simeq_{\beta'}&X'}
$$
such that $(\beta,\beta')\co\hat X\to X\times X'$ is also a fibration.
\end{lemma}

\begin{lemma}
\label{lemma-rigidify} Let
$$
\xymatrix{
&\hat A\ar[ld]_f\ar[d]^{\tilde f}\ar[rd]^{f'}\\
X&\ar[l]^\beta \hat X\ar[r]_{\beta'}&X' }
$$
be a homotopy commutative diagram in a closed model category. If
$\hat A$ is a cofibrant object, if $X$ and $X'$ are fibrant, and
if $(\beta,\beta')\co\hat X\to X\times X'$ is a fibration then
there exists a morphism $\hat f\co\hat A\to\hat X$ homotopic to
$\tilde f$ and making  the following diagram strictly commute
$$
\xymatrix{
&\hat A\ar[ld]_f\ar[d]^{\hat f}\ar[rd]^{f'}\\
X&\ar[l]^\beta \hat X\ar[r]_{\beta'}&X'. }
$$
\end{lemma}
\section{CDGA structures on mapping cones}\label{section-MC}
The aim of this section is to define a natural extension of the
$R$-DGmodule structure of some mapping cones to CDGA structures,
under certain dimension-connectivity hypotheses. 

\begin{defin}
\label{def-semitrivialCGA} Let $R$ be a CDGA and let $f\co X\to
R$ be a morphism of $R$-DGmodules. Consider the mapping cone
$C(f)=R\oplus_f sX$ and define a multiplication
$$
\mu\co C(f)\otimes C(f)\to C(f)
$$
by, for homogeneous elements $r,r'\in R$ and $x,x'\in X$,
\begin{itemize}
\item[(i)] $\mu(r\otimes r')=r.r'$
\item[(ii)] $\mu(r\otimes sx')=(-1)^{\deg(r)}s(r.x')$
\item[(iii)] $\mu(sx\otimes r')=(-1)^{\deg(x).\deg(r')}s(r'.x)$
\item[(iv)] $\mu(sx\otimes sx')=0$.
\end{itemize}
This multiplication defines a commutative graded algebra structure
(not necessarily differential) on $R\oplus_f sX$ that we call the
{\em semi-trivial CGA structure} on the mapping cone.
\end{defin}
This CGA structure on $C(f)$ is compatible with its $R$-module
structure in the sense that the module structure is induced by the
CGA map $R\hookrightarrow R\oplus_f sX$. It is important to notice
that in general  the multiplication $\mu$ defined above does not
define a CDGA structure on $C(f)$ because the Leibnitz rule on the
differential of the mapping cone is not necessarily satisfied.
However, we have the following lemmas.

\begin{lemma}\label{lemma-CDGAMC}%
Let $R$ be a CDGA and let $f\co X\to R$ be an $R$-DGmodule
morphism. Suppose that  $(sX)^{<k}=0$ and $(R\oplus sX)^{>2k}=0$ 
for some non negative integer $k$.
Then the mapping
cone $C(f)=R\oplus_f sX$ endowed with its semi-trivial
multiplication is a CDGA and the inclusion map $R\hookrightarrow
R\oplus_f sX$ is a CDGA-morphism.
\end{lemma}
\begin{proof}
This lemma is a special case of the next lemma with $I=0$ and $l=0$.
\end{proof}

\begin{lemma}\label{lemma-CDGAtruncMC}%
Let $R$ be a  CDGA, let $f\co X\to R$ be an $R$-DGmodule
morphism, and let $I\subset R\oplus_f sX$ be an $R$-DGsubmodule.
Suppose that 
$(sX)^{<k}=0$, $I^{\leq k-l}=0$, and $(R\oplus_f sX)^{\geq
2k-l+1}\subset I$
for non negative  
 integers $k$ and $l$. Then the semi-trivial multiplication $\mu$ on the
mapping cone $C(f)=R\oplus_fsX$ induces a multiplication on
$C(f)/I$ which endows this quotient with a CDGA-structure, and the
composition
$$
\xymatrix{R\,\ar@{^(->}[r]&R\oplus_fsX\ar[r]^{\pr}&C(f)/I}
$$
is a CDGA morphism.
\end{lemma}
\begin{proof}
We show first that $I$ is an ideal of the CGA $C(f)$ equipped with
its semi-trivial CGA structure. Since $I$ is an $R$-submodule of
$C(f)$ we have that $\mu(R\otimes I)=R.I\subset I$. On the other
hand, for degree reasons $\mu(sX\otimes I)\subset (C(f))^{\geq
2k-l+1}\subset I$. Therefore $\mu(C(f)\otimes I)\subset I$. Thus $I$
is a left ideal, hence a two-sided ideal because $\mu$ is graded
commutative.

This implies that the CGA structure on $C(f)$ induces a CGA
structure on the quotient $C(f)/I$. Denote by $\delta$ the
differential on the mapping cone $C(f)$ and by $\bar\delta$ the
induced differential on the quotient. To prove that
$(C(f)/I,\bar\delta)$ is a CDGA we have only to check the Leibnitz
formula. This will be a consequence of the following relation, for
$c,c'$ homogeneous elements in $R\oplus sX$:
\begin{equation}\label{equ-Leibniztrunc}%
\delta(\mu(c\otimes c'))-\mu(\delta(c)\otimes c')-(-1)^{|c|}
\mu(c\otimes\delta(c'))\,\in\,I.
\end{equation}
To prove \refequ{equ-Leibniztrunc} we study different cases. If
$c,c'\in R$ then the expression in \refequ{equ-Leibniztrunc} is
zero because $R$ is a DGA. If $c\in R$ and $c'\in sX$ then the
expression in \refequ{equ-Leibniztrunc} is zero because $\delta$
is a differential of $R$-DGmodule and the same is true if $c\in
sX$ and $c'\in R$ because $\mu$ is graded commutative. Finally if
$c,c'\in sX$ then the degree of the expression in
\refequ{equ-Leibniztrunc} is at least $2k+1\geq 2k-l+1$, therefore it belongs
to $I$.

This completes the proof that $C(f)/I$ is a CDGA. It is
straightforward to check that the map $R\to C(f)/I$ is a
CDGA-morphism.
\end{proof}

\begin{defin}\label{def-CDGAMC}
The CDGA-structures defined on the mapping cone $R\oplus_f sX$  in
Lemma \ref{lemma-CDGAMC} (respectively on the truncated mapping
cone $(R\oplus_f sX)/I$ in Lemma \ref{lemma-CDGAtruncMC}) is
called the \emph{semi-trivial CDGA structure}.
\end{defin}

Our last lemma gives a sufficient condition for some DGmodule map
between CDGA to be a CDGA morphism.
\begin{lemma}
\label{lemma-DGmodCDGAmap} Let $f\co A\to B$ be a
CDGA-morphism, let
$\xymatrix@1{A\quad\ar@{>->}[r]^-u&A\otimes\wedge X}$ be a relative
Sullivan algebra, and let $\hat f\co A\otimes\wedge X\to B$ be
an $A$-DGmodule morphism extending $f$. If $X^{<k}=0$ and
$B^{\geq2k}=0$ for some non negative integer $k$ then $\hat f$ is a CDGA morphism.
\end{lemma}
\begin{proof}
Since $A\otimes\wedge X$ and $B$ are graded commutative, $\hat f$
is a morphism of $A$-bimodules. The lemma follows from the fact
that for degree reasons $\hat f(A\otimes\wedge^{\geq 2}X)=0$.
\end{proof}

\section{\Topdegree or shriek map}
The aim of this section is to introduce the simple notion of a \emph{\topdegree map}
(which was called a \emph{shriek map} in early version of this paper).
A key result will be the existence and essential uniqueness of such \topdegree maps 
(Proposition \ref{prop-existsshriek}.)

We start with the definition and two examples.

\begin{defin}\label{def-shriek}%
Let $R$ be a DGA and assume that $H^*(R)$ is a connected
Poincar\'e duality algebra in dimension $n$. 
A \emph{\topdegree map of  $R$-DGmodule} 
 is an $R$-DGmodule map 
$\psi\co D\to R'$
such that $R'$ is weakly equivalent to $R$ and 
$H^n(\psi)$ is an isomorphism.
\end{defin}
\begin{example}
\label{ex-topdegreeshriek}
Suppose that $f\co V\hookrightarrow W$ is an embedding of
connected \emph{closed} oriented manifolds of codimension $k$.
Denote by [V] and [W] their homology orientation classes. We have
the classical cohomological shriek map (or Umkehr map, or Gysin
map, see \cite[VI.11.2]{Bredon})
$$
f^!\co s^{-k}H^*(V;\Bk)\to H^*(W;\Bk)
$$
characterized by the equation $f(s^{-k}v)\cap[W]=f_*(v\cap[V])$
(the $k$th-suspension is  here  only to make $f^!$ a degree
preserving
 map.) It is clear that $f^!$ is a map of $H^*(W)$-modules
and that it induces an isomorphism in
degree $n=\dim(W)$. Therefore $f^!$ is a \topdegree map of
$H^*(W)$-module (here
the differentials are supposed to be $0$).
\end{example}
\begin{example}
\label{ex-topdegreedual}
 Let $R$ be a DGA such that $H(R)$ is a connected
Poincar\'e duality  algebra in dimension $n$.  Let
$\phi\co R\to Q$ be a morphism of \emph{right} $R$-DGmodules
such that $H^0(\phi)$ is an isomorphism. Then $s^{-n}\#R$ is
quasi-isomorphic to $R$ and the map
$$s^{-n}\#\phi\co s^{-n}\#Q\to s^{-n}\#R$$
is a \topdegree map of (left) $R$-DGmodules.
\end{example}

To prove the existence and uniqueness of \topdegree maps we need first to study
further sets of homotopy classes of $R$-DGmodules.
For an integer $i$, denote by $\homset^i_R(P,N)$ the $\Bk$-module
of $R$-module maps of degree $i$ from $P$ to $N$ and set
$$\homset^*_R(P,N):=\oplus_{i\in\BZ} \homset^i_R(P,N).$$
We can define a degree $+1$ differential $\delta$ on this graded
$\Bk$-module by the formula $\delta(f) = d_N f-(-1)^{|f|}f d_P$.
The following  identification is well-known and we omit its proof
(e.g.\ \cite{FHT-gorenstein}):
\begin{lemma}\label{lemma-charsethmtpy}%
Let $R$ be a DGA, let $P$ be a cofibrant $R$-DGmodule, and let $N$
be an $R$-DGmodule. Then we have an isomorphism
$$[P,N]_R\cong H^0(\homset^*_R(P,N),\delta).$$
\end{lemma}

We have the following  important characterization of the set of
homotopy classes into a Poincar\'e duality algebra.
\begin{prop}
\label{prop-PDhmtpyclasses} Let $R$ be a DGA over a field $\Bk$
such that $H^*(R)$ is a connected Poincar\'e duality algebra in
 dimension $n$. Let $R'$ be an $R$-DGmodule weakly
equivalent to $R$ and let $P$ be a cofibrant $R$-DGmodule. Then
the map
$$
H^n\co [P,R']_R\to\homset_\Bk(H^n(P),H^n(R'))\,,\quad[f]\mapsto
H^n(f).
$$
is an isomorphism of $\Bk$-modules.
\end{prop}

\proof
Without any loss of generality we can suppose that $R'=R$ because
weak equivalences preserve each side of the isomorphism we want to
prove.

Since $H^n(R)\cong\Bk$ there exists a $\Bk$-DGmodule map
$
\epsilon_0\co R\to s^{-n}\Bk
$
inducing an isomorphism in $H^n$. Using the canonical isomorphism
$\# s^nR\cong s^{-n}\#R$ we can interpret $\epsilon_0$ as a
cocycle in $s^{-n}\#R$ and $[\epsilon_0]\not=0$ in
$H^0(s^{-n}\#R)\cong\#H^n(R)$.
 Since $R$ is also a \emph{right} $R$-DGmodule, we have a
structure of (left) $R$-DGmodule on $s^{-n}\#R$ (remember our
convention in Section \ref{section-toolkitApl}.) There is a unique
$R$-DGmodule map
$$\epsilon\co R\to s^{-n}\#R$$
sending $1\in R$ to $\epsilon_0$. Thus
$H^*(\epsilon)\co H^*(R)\to s^{-n}\#H^*(R)$
is an $H^*(R)$-module morphism which is an isomorphism in degree
$n$. By Poincar\'e duality of $H^*(R)$ this implies that
$H^*(\epsilon)$ is an isomorphism in every degree. Thus $\epsilon$
is a quasi-isomorphism.

Consider the adjunction isomorphism
$$\homset_R(P,\#R)\cong\homset_\Bk(P,\Bk)\,,\quad\phi\mapsto\hat\phi$$
where $\hat\phi\co P\to\Bk$ is defined by
$\hat\phi(x)=(\phi(x))(1)$ for $x\in P$ and $1$ the unit in $R$.
Combining this isomorphism with Lemma \ref{lemma-charsethmtpy}
 we get the following
sequence of isomorphisms
\begin{eqnarray*}
[P,R]_R&\cong& H^0(\homset_R(P,R))\\
&\stackrel{\epsilon_*}{\cong}&H^0(\homset_R(P,s^{-n}\#R))\\
&\cong&H^n(\homset_R(P,\#R))\\
&\cong&H^n(\homset_\Bk(P,\Bk))\\
&\cong&\homset_\Bk(H^n(P),s^n\Bk).
\end{eqnarray*}
Moreover it is straightforward to check that the following diagram
is commutative where the horizontal isomorphism is taken as the
previous sequence of isomorphisms:
$$
\xymatrix{[P,R]_R\ar[r]^-\cong\ar[rd]_{H^n}&\homset_\Bk(H^n(P),s^n\Bk)\\
&\homset_\Bk(H^n(P),H^n(R)).\ar[u]^\cong_{\epsilon^*_0}}
\eqno{\raise-38pt\hbox{\qed}}
$$

We establishes now the  existence and uniqueness (up to homotopy and a scalar multiple)
of \topdegree maps.
\begin{prop}\label{prop-existsshriek}
Let $R$ be a DGA such that $H^*(R)$ is a connected Poincar\'e
duality algebra in dimension $n$, let $R'$ be an
$R$-DGmodule weakly equivalent to $R$, and let $\hat D$ be a
cofibrant $R$-DGmodule such that $H^n(\hat D)\cong\Bk$. Then there
exists a \topdegree map of $R$-DGmodules
$$
\psi\co \hat D\to R'.
$$
Moreover if $\psi'\co \hat D\to R'$  is another \topdegree map
then there exists $u\in\Bk\takeaway\{0\}$
such that
$
[\psi]=u.[\psi']
$
in $[\hat D,R']_R$.
\end{prop}
\begin{proof}
By Proposition \ref{prop-PDhmtpyclasses} we have an isomorphism
$$
H^n\co[\hat D, R']_R\iso\homset_\Bk(H^n(\hat D),H^n(R')).
$$
Denote by $\isoset\left(H^n(\hat D),H^n(R')\right)$ the submodule
of $\homset_\Bk(H^n(\hat D),H^n(R'))$ consisting of isomorphisms.
Since $H^n(\hat D)\cong\Bk\cong H^n(R)$ there is an obvious
isomorphism
$$\isoset\left(H^n(\hat D),H^n(R')\right)\cong\Bk\takeaway\set{0}.
$$
Any  homotopy class $\psi\in[\hat D, R']_R$ corresponding to an
element of the non empty set $\isoset\left(H^n(\hat
D),H^n(R')\right)$ gives a \topdegree map, which proves the existence
part.

The uniqueness part is based on the same computation and left to the reader.
\end{proof}

We end this section by a lemma on sets of homotopy classes.
\begin{lemma}\label{lemma-stablehmtpyclasses}%
Let $A$ be a DGA, let $D$ be a cofibrant $A$-DGmodule, and let $X$
be an $A$-DGmodule. Suppose that there exist  integers $r\geq1$
and $m\geq0$ such that
\begin{itemize}
\item $H^{\leq r-1}(A)=H^0(A)=\Bk$, i.e.\ $A$ is cohomologically
$(r-1)$-connected,
\item $H^{<0}(X)=0$ and $H^{>m}(X)=0$, and
\item $H^{\leq m-r+1}(D)=0$.
\end{itemize}
 Then the map
$$H^*\co[D,X]_A\to \homset^0_\Bk(H^*(D),H^*(X))\,,\quad[f]\mapsto H^*(f)$$
is an isomorphism of $\Bk$-modules.

If moreover $r=1$ then
$[D,X]_A=0$.
\end{lemma}
\proof
We treat
separately the cases $r=1$ and $r\geq2$. Suppose first that $r=1$.
Then $H^{\leq m}(D)=0=H^{>m}(X)$. By standard obstruction theory
every $A$-DGmodule morphism $f\co D\to X$ is nullhomotopic.
Hence $[D,X]_A=0$. Moreover $\homset^0_\Bk(H^*(D),H^*(X))=0$ for
degree reasons. This proves the lemma for $r=1$.

Suppose that $r\geq2$. Using Lemma \ref{lemma-charsethmtpy} one
can prove that the $\Bk$-module $[D,X]_A$ remains unchanged if we
replace $D$, $X$, or $A$ by a cofibrant weakly
equivalent objects (see
\cite[Proposition A.4.(ii)]{FHT-gorenstein}.) Since $H^{\leq
1}(A)=\Bk$, we can replace the DGA $A$ by a minimal free model in
the sense of \cite[Appendix]{HalpLemcatDGA}, therefore we can
suppose that $A^{\leq r-1}=\Bk$. Next by replacing $D$ by a weakly
equivalent minimal semi-free
 $A$-DGmodule we can suppose that $D^{\leq m-r+1}=0$. Since $H^{>m}(X)=0$ and
$A$ is connected we can also assume  that $X^{>m}=0$.

Then, for degree reasons, the forgetful map
$
\phi^i\co\homset^i_A(D,X)\to \homset^i_\Bk(D,X)
$
is surjective for $i\geq -1$. Obviously $\phi^i$ is always
injective. Thus in  the following commutative diagram, the
horizontal maps are isomorphisms:
$$\xymatrix{
\homset_A^1(D,X)\ar[r]^{\cong}_{\phi^1}&\homset_\Bk^1(D,X)\\
\homset_A^0(D,X)\ar[r]^{\cong}_{\phi^0}\ar[u]_\delta&\homset_\Bk^0(D,X)\ar[u]_\delta\\
\homset_A^{-1}(D,X)\ar[r]^{\cong}_{\phi^{-1}}\ar[u]_\delta&\homset_\Bk^{-1}(D,X)\ar[u]_\delta.}
$$
This implies that $H^0(\phi)\co H^0(\homset^*_A(D,X),\delta)\to
H^0(\homset^*_\Bk(D,X),\delta)$ is an isomorphism. We conclude by
using Lemma \ref{lemma-charsethmtpy} and the obvious
identification
$$H^0(\homset^*_\Bk(D,X),\delta)\cong\homset^0_\Bk(H^*(D),H^*(X)).
\eqno{\qed}$$

\section{DGmodule model of a Poincar\'e embedding}\label{section-DGMod}

The aim of this section is to prove Theorem \ref{thm-DGmodnonconn}  and
Corollary \ref{corol-HWmodstruct}.

\begin{rmk}
Before proceeding with the proof of Theorem \ref{thm-DGmodnonconn}
we make a comment about the hypothesis on the model of a \emph{menorah}.
Indeed in that theorem we suppose that
$\set{\phi_k}_{1\leq k\leq c}$ is a model of the menorah
$\set{C^*(f_k)}_{1\leq k\leq c}$. As we pointed out in
Remark \ref{rmk-modelmenorah}, when $c\geq 2$ this is a stronger hypothesis
than asking for each $\phi_k$ to be a model of $C^*(f_k)$.
 We illustrate this fact by the
following example. Consider the torus $T=S^1\times S^1$ and denote
by $\dot{T}$ this torus with a small open disk removed, so that
$\dot{T}$ is a compact surface of genus $1$ with a circle for
boundary. Let $f\co S^1\hookrightarrow \dot{T}$ be an embedding
such that composed with the inclusion $\dot{T}\subset S^1\times
S^1$ it gives the inclusion of the first factor $S^1$ in
$S^1\times S^1$. Denote by $\dot{T_1}$ and $\dot{T_2}$ two copies
of $\dot{T}$ and let $f_k\co S^1\hookrightarrow\dot{T_k}$ be
the embeddings corresponding to $f$, $k=1,2$. Set
$W=\dot{T_1}\cup_{\del \dot{T}}\dot{T_2}$ which is a closed surface of
genus $2$. It is clear that the complement
$C:=W\takeaway(f_1(S^1)\amalg f_2(S^1))$ is connected. Consider
now the obvious automorphism $\phi$ of $W$ permuting $\dot{T_1}$
and $\dot{T_2}$. This automorphism is such that $\phi\circ
f_2=f_1$. By deforming slightly  $\phi$ into a diffeotopic
automorphism $\phi'$, we can suppose that $f'_2:=\phi'\circ
 f_2$ is an embedding of a circle closed
but disjoint from $f_1(S^1)$. Then $C':=W\takeaway(f_1(S^1)\amalg
f'_2(S^1))$ is not connected. Thus $C^*(C)$ and $C^*(C')$ do not
have the same DGmodule model since they have different
cohomologies. On the other hand $C^*(f'_2)$ and $C^*(f_2)$ do
admit the same model since they differ only by the automorphism
$C^*(\phi')$ of $C^*(W)$. The explanation of this apparent
contradiction is in the fact  that $\set{C^*(f_1),C^*(f_2)}$ and
$\set{C^*(f_1),C^*(f'_2)}$ do not admit a common model as
\emph{menorah} in the sense of Example \ref{examples-diagram}.
\end{rmk}

The proof of Theorem
 \ref{thm-DGmodnonconn}  consists of
a series of four lemmas.
Note first that by taking mapping cylinders we can assume
without loss of generality  that diagram \refequ{diag-mainsquare} 
of Definition \ref{def-Pemb} is a genuine push-out and that each map $i$, $k$, $f$, $l$ is a closed cofibration.
 
\begin{lemma}\label{lemma-C*MC}
With the same hypotheses as in Theorem \ref{thm-DGmodnonconn}
consider the inclusion map $\iota\co C^*(W,C)\to C^*(W)$. Then
the  commutative square $\BD'$  is weakly equivalent in
$C^*(W)$-DGMod to the following commutative square:

$\quad\quad\quad\quad\BD'':=\vcenter{\xymatrix@1{
C^*(\vrule width0pt depth6pt W)_{{}_{{}_{}}}\ar[r]^{f^*}\ar@{^(->}[d]&%
C^*(\vrule width0pt depth6pt P)_{{}_{{}_{}}}\ar@{^(->}[d]\\%
C^*(W)\oplus_{\iota}sC^*(W,C)\ar[r]^{f^*\oplus\id}&%
C^*(P)\oplus_{f^*\iota}sC^*(W,C)%
}}$
\end{lemma}
\begin{proof} Consider the following ladder of short exact sequences in
$C^*(W)$-DGMod
$$
\xymatrix{0\ar[r]&C^*(W,C)\ar[r]^{\iota}\ar[d]_{\simeq}^{f_0^*}&C^*(W)\ar[r]^{l^*}\ar[d]^{f^*}&
C^*(C)\ar[r]\ar[d]^{k^*}&0\\
0\ar[r]&C^*(P,\del T )\ar[r]^{\iota'}&C^*(P)\ar[r]^{i^*}& C^*(\del
T)\ar[r]&0.}
$$
By Mayer-Vietoris $f_0^*$ is a quasi-isomorphism and we have a weak
equivalence
$$ \id\oplus sf_0^*\co
\left(C^*(P)\oplus_{\iota'f_0^*}sC^*(W,C)\right)\quism
\left(C^*(P)\oplus_{\iota'}sC^*(P,\del T)\right).
$$
Thus in diagram $\BD''$ we can replace the right bottom DGmodule
by $C^*(P)\oplus_{\iota'}sC^*(P,\del T)$. To finish the proof apply the five lemma to deduce
that the map $k^*$ is weakly equivalent to the map induced between the mapping cones of $\iota$
and $\iota'$.
\end{proof}

Before stating the next two lemmas we need to introduce further
notation. Let $\del T_k$ be the union of the connected components of $\del T$ that are sent to $P_k$ by $i$. Set
 $C_k:=C\cup_{(\del T\takeaway\del T_k)}(P\takeaway P_k)$, which can 
 be interpreted as the complement of $P_k$ in $W$ since $W\simeq C_k\cup_{\del T_k}P_k$.
Define also the inclusion maps
$$
\iota_k\co C^*(W,C_k)\hookrightarrow C^*(W).
$$

In the next lemma we  build a convenient common model $\hat\phi_k$
of both $\phi_k$ and $f_k^*$.
\begin{lemma}\label{lemma-DGmodphihat}
With the hypotheses of Theorem \ref{thm-DGmodnonconn}  there
exists a cofibrant $A$-DGmodule $\hat R$, weak equivalences
$\alpha,\alpha'$, and, for each $k=1,\cdots,c$, an $A$-DGmodule
cofibration
$\xymatrix@1%
{\hat R\quad\ar@{>->}[r]^{\hat\phi_k}&\hat Q_k}%
$ and weak equivalences $\beta_k,\beta'_k$, making  the following
diagrams commute
$$
\xymatrix{%
R_{{}_{{}_{}}}\ar[d]_{\phi_k}& \hat
R_{{}_{{}_{}}}\ar[l]_{\alpha}^{\simeq}\ar[r]^{\alpha'}_{\simeq}\ar@{>->}[d]_{\hat\phi_k}&
C^*(W)_{{}_{{}_{}}}\ar[d]^{f^*_k}\\
Q_k& \hat Q_k\ar[l]^{\beta_k}_{\simeq}\ar[r]_{\beta'_k}^{\simeq}&
C^*(P_k), }
$$
and such that $(\alpha,\alpha')\co\hat R\to R\oplus C^*(W)$ and
$(\beta_k,\beta'_k)\co\hat Q_k\to Q_k\oplus C^*(P_k)$ are
surjective.
\end{lemma}
\begin{proof}
Let $\calS$ be the category shaping menorah's. 
Apply Lemma \ref{lemma-bisurj} in the category $A$-DGMod$^\calS$
to get a cofibrant menorah ${\set{\hat\phi_k}}_{1\leq k\leq c}$
and weak equivalences
$$
\xymatrix{%
\set{\phi_k}_{1\leq k\leq c}&
 {\quad\set{\hat\phi_k}}_{1\leq k\leq c}\quad
 \ar@{->>}[l]^-{\set{{(\alpha,\beta_k)}}_k}_{\simeq}
 \ar@{->>}[r]_-{\set{{(\alpha',\beta'_k)}}_k}^{\simeq}&
\set{f^*_k}_{1\leq k\leq c}}
$$
with the desired properties. In particular by the second part of
Proposition \ref{prop-CMdiagrams} the maps $\hat\phi_k\co \hat R\to \hat Q_k$
are cofibrations between cofibrant objects.
\end{proof}

\begin{lemma}\label{lemma-DGmodDhat}%
With the hypotheses of Theorem \ref{thm-DGmodnonconn} and with the
notation of Lemma \ref{lemma-DGmodphihat}, there exist for each
$k=1,\cdots,c$, a cofibrant $A$-DGmodule, $\hat D_k$, and weak
equivalences of $A$-DGmodules,
$$
\xymatrix{D_k&\ar[l]^{\simeq}_{\gamma_k}\hat
D_k\ar[r]^-{\gamma'_k}_-{\simeq}& C^*(W,C_k),}
$$
making the following diagram of isomorphisms commute
$$\xymatrix{%
H^n(D_k)\ar[d]_{H^n(\phishrk_k)}^\cong&\ar[l]_{H^n(\gamma_k)}^\cong
H^n(\hat
D_k)\ar[r]^{H^n(\gamma'_k)}_\cong&H^n(W,C_k)\ar[d]^{H^n(\iota_k)}_\cong\\
H^n(R)&\ar[l]^{H^n(\alpha)}_\cong H^n(\hat
R)\ar[r]_{H^n(\alpha')}^\cong& H^n(W).}
$$
\end{lemma}
\begin{proof}
Fix $k=1,\cdots,c$. By hypothesis $D_k$ is weakly equivalent as
an $A$-DGmodule to $s^{-n}\#C^*(P_k)$,
by Poincar\'e duality  to $C^*(P_k,\del T_k)$, and by Mayer-Vietoris
to $C^*(W,C_k)$. By Lemma \ref{lemma-bisurj}, we can find a
cofibrant $A$-DGmodule, $\hat D_k$, and weak equivalences of
$A$-DGmodules
$$
\xymatrix{D_k&\ar[l]^{\simeq}_{\gamma_k}\hat
D_k\ar[r]^-{\gamma''_k}_-{\simeq}& C^*(W,C_k).}
$$

By Lefschetz duality $H^n(W,C_k)\cong H_0(P_k)\cong\Bk$ and $H^n(\iota_k)$
is an isomorphism.
By  definition of a \topdegree map $H^n(\phishrk_k)$ is also an isomorphism.
Thus
the diagram appearing in the statement of the lemma, with
$\gamma''_k$ replacing $\gamma'_k$, is indeed a diagram of
isomorphisms. Since $H^n(\hat D_k)\cong H^n(\hat R)\cong\Bk$, the
two isomorphisms
$$
H^n(\alpha)^{-1}H^n(\phishrk_k)H^n(\gamma_k)
\textrm{\,\,\,\,and\,\,\,\,}
H^n(\alpha')^{-1}H^n(\iota_k)H^n(\gamma''_k)
$$
differ only by a multiplicative constant
$u\in\Bk\takeaway\set{0}$. Set $\gamma'_k:=u.\gamma''_k$ which is
also a weak equivalence of $A$-DGmodules. Then the diagram of
isomorphisms of the statement commutes.
\end{proof}

Recall the notion of  model of a composite from Example
\ref{examples-diagram}.

\begin{lemma}
\label{lemma-DGmodphiphishriek} With the hypotheses of Theorem
\ref{thm-DGmodnonconn}, the composite
$$
\xymatrix{%
D\ar[r]^{\phishrk}&R\ar[r]^\phi&Q}
$$
is an $A$-DGmodule model of the composite
$$
\xymatrix{%
C^*(W,C)\ar[r]^{\iota}&C^*(W)\ar[r]^{f^*}&C^*(P).}
$$
\end{lemma}
\begin{proof}
Consider all the morphisms and DGmodules built in Lemma
\ref{lemma-DGmodphihat} and Lemma \ref{lemma-DGmodDhat}. Fix
$k=1,\cdots,c$. Take a
lifting of $A$-DGmodules $\hat\phishrk_k\co\hat D_k\to\hat R$ of
$\phishrkk\gamma_k$ along the acyclic fibration $\alpha$, so that
\begin{equation}
\label{equ-alphaphi1} \alpha\hat\phishrkk=\phishrk_k\gamma_k
\end{equation}
which is a \topdegree map. Also $\alpha'\hat\phishrk_k$ and $\iota_k\gamma'_k$ are
\topdegree maps with values in $C^*(W)$. By Proposition \ref{prop-existsshriek}
there are homotopic up to a multiplicative scalar $u\not=0$ and 
 Lemma \ref{lemma-DGmodDhat} and \refequ{equ-alphaphi1} imply that $u=1$. Thus
\begin{equation}\label{equ-alphaphi2}%
\alpha'\hat\phishrk_k\simeq_A\iota_k\gamma'_k.
\end{equation}

Set  $\hat D:=\oplus_{k=1}^c\hat D_k$,
$\gamma:=\oplus_{k=1}^c\gamma_k$,
$\gamma':=\oplus_{k=1}^c\gamma'_k$, and
$\hat\phishrk:=\sum_{k=1}^c\hat\phishrk_k$. Since the $P_k$'s are
pairwise disjoint, we have an identification
$C^*(W,C)=\oplus_{k=1}^c C^*(W,C_k)$. Equations
$(\ref{equ-alphaphi1})$ and $(\ref{equ-alphaphi2})$ yield to the
following homotopy commutative diagram in $A$-DGMod
$$
\xymatrix{%
D\ar[d]_{\phishrk}& \hat
D\ar[l]_{\gamma}^{\simeq}\ar[r]^-{\gamma'}_-{\simeq}\ar[d]_{\hat\phishrk}&
C^*(W,C)\ar[d]^{\iota}\\
R& \hat R\ar[l]_{\alpha}^{\simeq}\ar[r]^-{\alpha'}_-{\simeq}&
C^*(W). }
$$
Since $(\alpha,\alpha')$ is a fibration we can suppose by Lemma
\ref{lemma-rigidify} that $\hat\phishrk$ has been chosen such that
the above diagram is strictly commutative. Gluing this diagram
with that built in Lemma \ref{lemma-DGmodphihat} we get a
commutative diagram of $A$-DGmodules
$$
\xymatrix{%
D\ar[d]_{\phishrk}& \hat
D\ar[l]_{\gamma}^{\simeq}\ar[r]^-{\gamma'}_-{\simeq}\ar[d]_{\hat\phishrk}&
C^*(W,C)\ar[d]^{\iota}\\
R\ar[d]_{\phi}& \hat
R\ar[l]_{\alpha}^{\simeq}\ar[r]^-{\alpha'}_-{\simeq}\ar[d]_{\hat\phi}&
C^*(W)\ar[d]^{f^*}\\
Q& \hat Q\ar[l]^{\beta}_{\simeq}\ar[r]_-{\beta'}^-{\simeq}& C^*(P)
}
$$
and the lemma is proved.\end{proof}

Collecting the four previous lemmas we achieve the proof of
Theorem \ref{thm-DGmodnonconn} and its corollary.

\begin{proof}[Proof of Theorem \ref{thm-DGmodnonconn}]
Recall  diagrams $\BD$ and $\BD'$ defined in Theorem
\ref{thm-DGmodnonconn} and diagram $\BD''$ defined in Lemma
\ref{lemma-C*MC}. Using Lemma \ref{lemma-DGmodphiphishriek} and 
taking mapping cones we deduce that the diagrams $\BD$ and
$\BD''$ are weakly equivalent in $A$-DGMod. By  Lemma
\ref{lemma-C*MC} diagrams $\BD'$ and $\BD''$ are also weakly
equivalent in $A$-DGMod. 
\end{proof}
\begin{proof}[Proof of Corollary \ref{corol-HWmodstruct}]
By Example \ref{ex-topdegreedual}  $s^{-n}\#\phi$
is a \topdegree map of right $A$-DGmodules.
Theorem \ref{thm-DGmodnonconn} implies that 
$s^{-n}\#R\oplus_{s^{-n}\#\phi}ss^{-n}\#Q$ is a right $A$-DGmodule model
of $C^*(C)$. Therefore their homologies are isomorphic as right $H^*(W)$-modules
and by commutativity also as left modules.

Since $H^{>m}(P)=0$ and by Lefschetz duality, $H_{<n-m}(W,C)=0$ and $H^i(W)\to H^i(C)$ is
an isomorphism for $i<n-m-1$. Therefore if $x.y$ is
a product in $H^*(C)$ that is not determined by the $H^*(W)$-module structure then
$\deg(x),\deg(y)\geq n-m-1$. Hence $\deg(x.y)\geq 2(n-m-1)\geq n-r$. Since $H^*(f)$ is
$r$-connected we have that $H_{\leq r}(W,P)=0$ and by Lefschetz duality $H^{\geq n-r}(C)=0$.
Therefore $x.y=0$. 
\end{proof}
\section{CDGA model of a Poincar\'e embedding in the stable case}
\label{section-stableCDGA}  
In this section we give a proof of  Theorem \ref{thm-stableCDGA}.
Here is an overview of that proof.
\begin{enumerate}
\item We want to show that the diagrams $\BD$ and $\BD'$ are weakly equivalent as commutative squares of
CDGA. By Theorem \ref{thm-DGmodnonconn} we already know that they are
weakly equivalent in a certain category of DGmodules.
\item
We will build a convenient  common CDGA model $\xymatrix@1{\hat
R\quad\ar@{>->}[r]^{\hat\phi}&\hat Q}$   of both $\phi\co R\to
Q$ and  $f^*\co\Apl(W)\to\Apl(P)$. We can then consider the
category of ``$\hat\phi$-DGmodules'' whose objects consist of maps
of $\hat R$-DGmodules $M\to N$ such that $N$ is also equipped with
a $\hat Q$-DGmodule compatible with its $\hat R$-DGmodule
structure through the map $\hat\phi$. The morphisms of this
category consist of certain commutative squares that we call
\emph{$\hat\phi$-squares} (see Definition
\ref{def-phihatsq}). In particular the diagrams $\BD$ and $\BD'$
will be $\hat\phi$-squares.
\item A refinement of the arguments of Theorem \ref{thm-DGmodnonconn}
will show that the diagrams $\BD$ and $\BD'$ are weakly equivalent
not only as squares in the category of $\hat R$-DGmodules but also
as $\hat\phi$-squares, which means that the weak equivalences
between the right sides of diagrams $\BD$ and $\BD'$ will be of
$\hat Q$-DGmodules (Lemma \ref{lemma-thetatheta''}.)
\item Using the results of Section \ref{section-MC}
(notably  Lemma \ref{lemma-DGmodCDGAmap}), we will
show that this weak equivalence of $\hat\phi$-squares between
$\BD$ and $\BD'$ is indeed a weak equivalence of CDGA squares.
\end{enumerate}

Let's move to the details by establishing a series of lemmas.
Note first that without loss of generality we can assume that \refequ{diag-mainsquare}
is a genuine push-out and that $f$ induces a map of pairs
$f_0\co(P,\del T)\to (W,C)$.
In the next lemma we build a common model $\xymatrix@1{\hat
R\quad\ar@{>->}[r]^{\hat\phi}&\hat Q}$ of both $\phi$ and $f^*$.
\begin{lemma}\label{lemma-CDGAphihat}
With the hypotheses of Theorem \ref{thm-stableCDGA}  there exists
a cofibrant CDGA  $\hat R$, a relative Sullivan algebra
$\xymatrix@1{\hat R\quad\ar@{>->}[r]^{\hat\phi}&\hat Q}$, and a
commutative diagram of CDGA where horizontal arrows are weak
equivalences
$$
\xymatrix{%
R\ar[d]_{\phi}& \hat
R\ar[l]_{\alpha}^{\simeq}\ar[r]^-{\alpha'}_-{\simeq}\ar[d]_{\hat\phi}&
\Apl(W)\ar[d]^{f^*}\\
Q& \hat Q\ar[l]^{\beta}_{\simeq}\ar[r]_-{\beta'}^-{\simeq}&
\Apl(P), }
$$
and  $(\alpha,\alpha')\co\hat R\to R\oplus \Apl(W)$
and $(\beta,\beta')\co\hat Q\to Q\oplus \Apl(P)$ are
surjections.
\end{lemma}
\begin{proof}
This is a consequence of Lemma \ref{lemma-bisurj} in the
the category of maps in CDGA, of the second part of Proposition
\ref{prop-CMdiagrams}, and of the fact that every CDGA cofibration is a retract of a
Sullivan relative algebra. Alternatively the lemma can be proved using standard
techniques of \cite{FHT-RHT}. 
 \end{proof}

Our next lemma gives a replacement $\bar R$ of $\hat R$ that fibres on
different DGmodules.
\begin{lemma}\label{lemma-Rbar}
With the hypotheses of Theorem \ref{thm-stableCDGA} and the
notation of Lemma \ref{lemma-CDGAphihat} there exists a
factorization of $\hat R$-DGmodules of $(\alpha,\alpha',\hat\phi)$
into an acyclic cofibration $\rho$ followed by a fibration
$(\bar\alpha,\bar\alpha',\bar\phi)$ as follows:
$$
\xymatrix{ \hat
R\,\,\ar@{>->}[r]_{\simeq}^\rho\ar[d]_{(\alpha,\alpha',\hat\phi)}&
\bar
R\ar@{->>}[ld]^{(\bar\alpha,\bar\alpha',\bar\phi)}\\
R\oplus\Apl(W)\oplus\hat Q&}
$$
\end{lemma}
\begin{proof}
The existence of such a factorization is one of the axioms of the
closed model structure on the category $\hat R$-DGMod.
\end{proof}

In the following lemma we give  a common model $\hat D$ of both
$D$ and $\Apl(P,\del T)$.
\begin{lemma}\label{lemma-stableDhat}
With the hypotheses of Theorem \ref{thm-stableCDGA} and with the
notation of Lemmas \ref{lemma-CDGAphihat} and \ref{lemma-Rbar},
there exists a cofibrant $\hat Q$-DGmodule, $\hat D$, and weak
equivalences of $\hat Q$-DGmodules,
$$
\xymatrix{D&\ar[l]^{\simeq}_{\gamma}\hat
D\ar[r]^-{\gamma'}_-{\simeq}& \Apl(P,\del T),}
$$
making the following diagram of isomorphisms commute
$$\xymatrix{%
H^n(D)\ar[d]_{H^n(\phishrk)}^\cong &\ar[l]_{H^n(\gamma)}^\cong
H^n(\hat
D)\ar[r]^{H^n(\gamma')}_\cong &H^n(P,\del T)&H^n(W,C)\ar[l]_{H^n(f_0)}^\cong \ar[ld]^{H^n(\iota)}_\cong \\
H^n(R)&\ar[l]^{H^n(\bar\alpha)}_\cong H^n(\bar
R)\ar[r]_{H^n(\bar\alpha')}^\cong & H^n(W).}
$$
Moreover $\hat D$ is also a cofibrant $\hat R$-DGmodule and there
exists an $\hat R$-DGmodule weak equivalence
$$
\gamma''\co\hat D\quism\Apl(W,C)
$$
making  the following diagram commute
$$
\xymatrix{\hat
D\ar[r]^-{\gamma''}_-\simeq\ar[rd]_{\gamma'}^\simeq&\Apl(W,C)\ar[d]^{f^*_0}_\simeq\\
&\Apl(P,\del T).}
$$
\end{lemma}
\begin{proof}
The  proof of the first part of the lemma
is similar to the proof of Lemma \ref{lemma-DGmodDhat}.

For the second part of the lemma, note that by \cite[Lemma
14.1]{FHT-RHT} $\hat D$ is a cofibrant $\hat R$-DGmodule because
it is a cofibrant $\hat Q$-DGmodule and because
$\hat\phi\co\hat R\to\hat Q$ is a relative Sullivan algebra.
Also $f^*_0$ is a surjective quasi-isomorphism. We take
 $\gamma''$ as a lift of $\gamma'$ along the acyclic fibration $f^*_0$.
\end{proof}

\begin{lemma}\label{lemma-CDGAphitildeshriek}
With  the hypotheses of Theorem \ref{thm-stableCDGA} and with the
notation of Lemmas \ref{lemma-CDGAphihat}--\ref{lemma-stableDhat},
there exists an $\hat R$-DGmodule morphism
$$
\tilde\phishrk\co\hat D\to\bar R
$$
making  the following diagram homotopy commute in $\hat R$-DGMod
$$\xymatrix{
D\ar[d]_{\phishrk}
&\ar[l]_{\gamma}\hat
D\ar[r]^-{\gamma''}\ar[d]_{\tilde\phishrk}
&
\Apl(W,C)\ar[d]^{\iota}\\
R&\ar[l]^{\bar\alpha}\bar R\ar[r]_-{\bar\alpha'}&\Apl(W).}$$
\end{lemma}
\begin{proof}
The argument is the same as in the beginning of the
proof of Lemma \ref{lemma-DGmodphiphishriek}.
\end{proof}

We build now a $\hat Q$-DGmodule common model $\chi$ both of
$\phi\phishrk=0$ and of $\iota'\co \Apl(P,\del T)\to\Apl(P)$.
\begin{lemma}\label{lemma-chi}
With  the hypotheses of Theorem \ref{thm-stableCDGA} and with the
notation of Lemmas
\ref{lemma-CDGAphihat}--\ref{lemma-CDGAphitildeshriek}, the
composite $\phi\phishrk$ is a $Q$-DGmodule morphism and there exists
a $\hat Q$-DGmodule morphism
$$
\chi\co\hat D\to\hat Q
$$
making  the following diagram commute in $\hat Q$-DGMod
$$
\xymatrix{ D\ar[d]_{\phi\phishrk}&\ar[l]_\gamma\hat
D\ar[d]_{\chi}\ar[r]^-{\gamma'}&\Apl(P,\del T)\ar[d]^{\iota'} \\
Q&\ar[l]^{\beta}\hat Q\ar[r]_-{\beta'}&\Apl(P).}
$$
Moreover $\chi\simeq_{\hat R}\bar\phi\tilde\phishrk$.
\end{lemma}
\begin{proof}
Notice that for degree reasons $\phi\phishrk=0$, therefore it is a
morphism of $Q$-DGmodules. Applying Lemma \ref{lemma-stablehmtpyclasses} with $r=1$
we get that $[\hat D,Q]_{\hat Q}=0=[\hat D,\Apl(T)]_{\hat Q}$.
Therefore the diagram of the statement with $0$ replacing $\chi$
is homotopy commutative in $\hat
Q$-DGMod.
Since $(\beta,\beta')$ is a fibration, Lemma \ref{lemma-rigidify}
permits to replace the zero map by a
homotopic $\hat Q$-DGmodule morphism $\chi$
making the diagram strictly commute.

We have also by Lemma \ref{lemma-stablehmtpyclasses} that $ [\hat
D,\hat Q]_{\hat R}=0$, hence $\chi\simeq_{\hat
R}\bar\phi\tilde\phishrk$.
\end{proof}
\begin{lemma}\label{lemma-factorchi}
With  the hypotheses of Theorem \ref{thm-stableCDGA} and with the
notation of Lemmas \ref{lemma-CDGAphihat}--\ref{lemma-chi}, there
exists a morphism of $\hat R$-DGmodule
$
\bar\phishrk\co\hat D\to\bar R
$
making  both of the following diagrams commute
$$\xymatrix{
D\ar[d]_{\phishrk}&\ar[l]_{\gamma}\hat
D\ar[r]^{\gamma''}\ar[d]_{\bar\phishrk}&
\Apl(W,C)\ar[d]^{\iota}&\textrm{\,\,\,\,and\,\,\,\,}
&\hat D\ar[r]^{\bar\phishrk}\ar[rd]_\chi&\bar R\ar[d]^{\bar\phi}\\
R&\ar[l]^{\bar\alpha}\bar R\ar[r]_{\bar\alpha'}&\Apl(W)&&&\hat Q.}
$$
\end{lemma}
\begin{proof}
By Lemmas  \ref{lemma-CDGAphitildeshriek} and \ref{lemma-chi}, we
have the following homotopy commutative diagram in $\hat R$-DGMod
$$\xymatrix{&\bar R\ar@{->>}[d]^{(\bar\alpha,\bar\alpha',\bar\phi)}\\
\hat
D\ar[ur]^{\tilde\phishrk}\ar[r]_-{(\phishrk\gamma,\iota\gamma'',\chi)}&\quad
R\oplus\Apl(W)\oplus\hat Q. }
$$
Since $(\bar\alpha,\bar\alpha',\bar\phi)$ is a fibration and $\hat
D$ is cofibrant, a standard argument in closed model
categories shows that we can replace $\tilde\phishrk$ by a homotopic
map $\bar\phishrk$ making the diagram strictly commute.
\end{proof}

As we have explained in the overview of the proof, in order to
prove that  diagrams $\BD$ and $\BD'$ of Theorem
\ref{thm-stableCDGA} are weakly equivalent in CDGA, we will first
prove that there are weakly equivalent as ``$\hat\phi$-squares''
that we define now. To give a meaning to this assertion we
could define a genuine closed model structure on the category of
$\hat\phi$-squares. Instead of doing so we prefer to introduce the
following \emph{ad hoc} definition of weakly equivalent
$\hat\phi$-squares.
\begin{defin}\label{def-phihatsq}
Let $\hat\phi\co\hat R\to\hat Q$ be a CDGA morphism.
\begin{enumerate}
\item[(i)] By a \emph{$\hat\phi$-square} we mean a commutative square of
$\hat R$-DGmodules
$$
\xymatrix{M\ar[r]^\psi\ar[d]_f&N\ar[d]^g\\
M'\ar[r]_{\psi'}&N'}
$$
such that $N$ and $N'$ have also a structure of $\hat Q$-DGmodule
compatible with their $\hat R$-DGmodule structure through
$\hat\phi$ and such that the right map $g$ is a $\hat
Q$-DGmorphism.
\item[(ii)] A \emph{morphism of $\hat\phi$-squares} is a morphism,
$\Theta$,  of commutative squares in $\hat R$-DGmodules between
two
$\hat\phi$-squares
$$\vcenter{\xymatrix@1{%
M\ar[r]^\psi\ar[d]_f &%
N\ar[d]^g\\
M'\ar[r]_{\psi'}&%
N'}}%
\quad\stackrel{\Theta}{\to}\quad%
\vcenter{\xymatrix@1{%
X\ar[r]^\omega\ar[d]_p&%
Y\ar[d]^q\\
X'\ar[r]_{\omega'}&%
Y'}}$$
of the form
$\Theta={\left(\begin{array}{cc}\mu&\nu\\\mu'&\nu'\end{array}\right)}$
 such that $\nu$ and $\nu'$ are also morphisms of $\hat
Q$-DGmodules.

\item[(iii)] A morphism  $\Theta$ of $\hat\phi$-squares is called a
\emph{fibration} (resp.\ \emph{a weak equivalence}) if each of the
morphisms $\mu,\mu',\nu,\nu'$ is a surjection (resp.\
quasi-isomorphism). A morphism of $\hat\phi$-squares which is
both a fibration and a weak equivalence is called an \emph{acyclic
fibration}.
\end{enumerate}
\end{defin}

Recall the diagrams $\BD$ and $\BD'$ from the statement of Theorem
\ref{thm-stableCDGA}.
\begin{lemma}\label{lemma-phihatsq}
With  the hypotheses of Theorem \ref{thm-stableCDGA} and with the
notation of Lemmas \ref{lemma-CDGAphihat}--\ref{lemma-factorchi},
Diagrams $\BD$ and $\BD'$ are both commutative squares in CDGA and
$\hat\phi$-squares. The following  diagram
$$
\bar \BD:=%
\vcenter{\xymatrix{%
\bar R_{{}_{{}_{}}}\ar[r]^{\bar\phi}\ar@{^(->}[d]&%
\hat Q_{{}_{{}_{}}}\ar@{^(->}[d]\\%
\bar R\oplus_{\bar\phishrk}s\hat D\ar[r]_{\bar\phi\oplus\id}&%
\hat Q\oplus_\chi s\hat D%
}}%
$$
is a $\hat\phi$-square.
\end{lemma}
\begin{proof}
The CDGA structure on the mapping cones of the bottom side of
Diagram $\BD$ are the semi-trivial CDGA structures, which exist by
Lemma \ref{lemma-CDGAMC}. From this it is clear that $\BD$ is a
commutative square of CDGA, as well as $\BD'$. They are also
$\hat\phi$-squares with $\hat R$- and $\hat Q$-DGmodule structures
induced by the maps $\alpha$, $\alpha'$, $\beta$, and $\beta'$.

Using the fact that $\chi$ is a $\hat Q$-DGmodule morphism it is
immediate to check that $\bar \BD$ is a $\hat\phi$-square.
\end{proof}
\begin{lemma}\label{lemma-thetatheta''}
With  the hypotheses of Theorem \ref{thm-stableCDGA} and with the
notation of Lemmas \ref{lemma-CDGAphihat}--\ref{lemma-phihatsq},
there exist acyclic fibrations of $\hat\phi$-squares
$$
\xymatrix{\BD&%
\ar@{->>}[l]_{\Theta}^{\simeq}\bar\BD\ar@{->>}[r]^{\Theta'}_\simeq&%
\BD'.}
$$
\end{lemma}
\begin{proof}
Using the different maps constructed in our previous series of
lemmas we will describe  these two acyclic fibrations
explicitly.
Consider the following commutative square
$$\BD'':=%
\vcenter{\xymatrix@1{%
\Apl(\vrule width0pt depth6pt W)_{{}_{{}_{}}}\ar[r]^{f^*}\ar@{^(->}[d]&%
\Apl(\vrule width0pt depth6pt P)_{{}_{{}_{}}}\ar@{^(->}[d]\\%
\Apl(W)\oplus_\iota s\Apl(W,C)\ar[r]_{f^*\oplus sf^*_0}%
&\Apl(P)\oplus_{\iota'} s\Apl(P,\del T).}}$$
Using the fact that $\iota'\co \Apl(P,\del T)\to \Apl(P)$ is a
morphism of $\Apl(P)$-DGmodules, hence of $\hat Q$-DGmodules, we
see that $\BD''$ is a diagram of $\hat\phi$-squares. Clearly
$\Theta''':=%
\left(%
\begin{array}{cc}%
\id&\id\\%
l^*\oplus 0&i^*\oplus 0%
\end{array}%
\right)\co\BD''\to\BD'%
$ is a surjection, and an argument analogous to that of Lemma
\ref{lemma-C*MC} shows that it is a weak equivalence. Hence
$\Theta'''$ is an acyclic fibration. We have  another acyclic
fibration
$\xymatrix@1{\Theta''\co\bar\BD\ar@{->>}[r]^\simeq&\BD''}$
given by
$\Theta'':=%
\left(%
\begin{array}{cc}%
\bar\alpha'&\beta'\\%
\bar\alpha'\oplus s\gamma''&\beta'\oplus s\gamma'%
\end{array}%
\right).%
$ Then $\Theta':=\Theta'''\Theta''$ is one of the required acyclic
fibration. The other one is given by $
\Theta:=%
\left(%
\begin{array}{cc}%
\bar\alpha&\beta\\%
\bar\alpha\oplus s\gamma&\beta\oplus s\gamma%
\end{array}%
\right).%
$
\end{proof}
We sketch now an overview of the end of the proof of the Theorem.
In the next lemma we build an intermediate commutative square,
$\hat\BD$, which is a CDGA model of $\BD'$. Moreover $\hat\BD$ is
also a ``cofibrant $\hat\phi$-square'', therefore by lifting
along the quasi-isomorphisms $\Theta$ and $\Theta'$ we will deduce
that $\hat \BD$ is a model of $\hat\phi$-square of $\BD$.
 Finally a degree argument will imply that
this $\hat\phi$-square quasi-isomorphism
$\hat\Theta\co\hat\BD\simeq\BD$ is in fact of CDGA and this
will prove that $\BD$ and $\BD'$ are weakly equivalent CDGA
squares. Let's move  to the details.

\begin{lemma}\label{lemma-thetahat''}
With  the hypotheses of Theorem \ref{thm-stableCDGA} and with the
notation of Lemmas \ref{lemma-CDGAphihat} and
\ref{lemma-phihatsq}, there exists a commutative square in CDGA,
$$\hat\BD:=%
\vcenter{\xymatrix@1{%
\hat {R}_{{}_{{}_{{}_{{}_{}}}}}\quad\ar@{>->}[r]^{\hat\phi}\ar@{>->}[d]_u&%
\hat Q_{{}_{{}_{{}_{{}_{}}}}}\ar@{>->}[d]^v\\
\hat R\otimes\wedge X\quad\ar@{>->}[r]_-{\hat\psi}&%
\hat Q\otimes \wedge X\otimes\wedge Y}} $$
where $\hat\phi$, $\hat\psi$, $u$, and $v$  are cofibrations, 
together with a weak equivalence both of CDGA-squares and of
$\hat\phi$-squares
$\hat\Theta'\co\hat\BD\quism \BD'$.
Moreover $X$ and $Y$ can be chosen such that such that
$X^{<n-m-1}=Y^{<n-m-2}=0$. If $H^1(f;\BQ)$ is injective
we can also assume that $Y^{<n-m-1}=0$.
\end{lemma}

\begin{proof}
By taking a {minimal} relative Sullivan algebra of $l^*\alpha'$
we get a commutative diagram of CDGA
$$\xymatrix{
\hat Q\ar[d]_{\beta'}&\ar[l]_{\hat\phi}\hat R\quad\ar[d]_{\alpha
'}\ar@{>->}[r]^u&\hat R\otimes \wedge X\ar[d]^{\lambda'}\\
\Apl(P)&\ar[l]^{f^*}\Apl(W)\ar[r]_{l^*}&\Apl(C).}
$$
Consider the push-out $\hat Q\otimes\wedge X$ of the top line of
the above diagram. By the universal property of the push-out, this
diagram induces a CDGA map
$$\bar\mu'_0\co\hat Q\otimes\wedge X\to\Apl(\del T).$$ The
latter map can be factored into a minimal relative Sullivan
algebra followed by a quasi-isomorphism,
$$\xymatrix{\hat Q\otimes\wedge X\,\,\ar@{>->}[r]^-v&
\hat Q\otimes\wedge X\otimes\wedge Y\ar[r]^-{\mu'}_\simeq&\Apl(\del
T).}
$$
It is immediate to check that the matrix
$\hat\Theta'=\left(\begin{array}{cc}
\alpha'&\beta'\\\lambda'&\mu'\end{array}\right)$ is a weak
equivalence of CDGA-squares and of $\hat\phi$-squares.

We prove now that $X^{<n-m-1}=0$. Since $i$ is $(n-m-1)$-connected
a Mayer-Vietoris argument implies that $H_*(l)$ is $(n-m-1)$-connected. Therefore the same is
true for the map $u$ and by minimality we get that $X^{<n-m-1}=0$.

The model $\hat Q\to \hat Q\otimes\wedge X\otimes\wedge Y$ of $i^*$
is cohomologically $(n-m-1)$-connected and since $X^{<n-m-1}=0$, minimality
implies that $Y^{<n-m-2}=0$.

Assume that $H^1(f)$ is injective. Thus $H^*(f)$ is $1$-connected and $H^*(l)$ is
 $(n-m-1)$-connected . By a rational Blackers-Massey argument we deduce that
 $\hat Q\otimes\wedge X\to \hat Q\otimes\wedge X\otimes\wedge Y$
 is cohomologically $(n-m-1)$-connected. By minimality
we get that $Y^{<n-m-1}=0$.
\end{proof}
\begin{lemma}\label{lemma-bartheta}
With  the hypotheses of Theorem \ref{thm-stableCDGA} and with the
notation of Lemma \ref{lemma-thetahat''}, the commutative squares of CDGA
$\hat\BD$ and $\BD$ are weakly equivalent.
\end{lemma}

\begin{proof}
The proof is in two steps:

(i)\qua We will show that the morphism $\hat\Theta'$ constructed in
Lemma \ref{lemma-thetahat''} lifts along the acyclic fibration
$\Theta'$ of Lemma \ref{lemma-thetatheta''} to a weak equivalence
of $\hat\phi$-squares
$\bar\Theta\co\hat\BD\quism\bar\BD$.

(ii)\qua using the map $\Theta$ constructed in Lemma
\ref{lemma-thetatheta''} we will show that the composite
$\hat\Theta:=\Theta\bar\Theta$ is a quasi-isomorphism of CDGA.

(i)\qua The lift will be of the form
$\bar\Theta:=\left(\begin{array}{cc}%
\rho&\id\\%
\bar\lambda&\bar\mu%
\end{array}\right)$, where $\rho$ was defined in Lemma \ref{lemma-Rbar}.
We need only to build the maps $\bar\lambda$ and $\bar\mu$, and
for this we will use the maps $\lambda'$ and $\mu'$ constructed in
the proof of Lemma \ref{lemma-thetahat''}. We have the following
solid commutative  diagram of $\hat R$-DGmodules
$$
\xymatrix{ \vrule width0pt depth6pt \hat R\ar@{>->}[d]_u\ar[r]^\rho&\bar R\ar@{^(->}[r]&
\bar R\oplus_{\bar\phishrk}s\hat
D\ar@{->>}[d]_{\simeq}^{l^*\bar\alpha'\oplus 0}\\
\hat R\otimes\wedge
X\ar@{-->}[rru]^{\bar\lambda}\ar[rr]^{\simeq}_{\lambda'}&&\Apl(C).}
$$
Since $u$ is a relative Sullivan algebra, it is an $\hat
R$-DGmodule cofibration. Then,  $l^*\bar\alpha\oplus 0$ being an
acyclic fibration, there exists a lift of $\hat R$-DGmodules,
$\bar\lambda$, making both triangles of the diagram commute.

We can define a map of $\hat Q$-DGmodules
$$
\bar\mu_0\co\hat Q\otimes\wedge X\to\hat Q\oplus_\chi s\hat D
$$
by the formula
$$\bar\mu_0(q\otimes\omega)=q.(\bar\phi\oplus\id)(\bar\lambda(1\otimes\omega)),$$
for $q\in\hat Q$ and $\omega\in\wedge X$. It is immediate to check
that $\bar\mu_0$ is a $\hat Q$-DGmodule morphism and that the
following solid diagram of $\hat Q$-DGmodules commutes:
$$\xymatrix{
\hat Q\otimes\wedge X\ar@{>->}[d]\ar[r]^{\bar\mu_0}&\hat
Q\oplus_\chi
s\hat D\ar@{->>}[d]_\simeq^{i^*\beta'\oplus 0}\\
\hat Q\otimes\wedge X\otimes\wedge
Y\ar@{-->}[ru]^{\bar\mu}\ar[r]_{\mu'}^\simeq&\Apl(\del T).}
$$
Therefore there exists a lift, $\bar\mu$, of $\hat Q$-DGmodules
making  both triangles of the diagram commute.

It is immediate to check that $\bar\Theta:=\left(\begin{array}{cc}
\rho&\id\\\bar\lambda&\bar\mu\end{array}\right)$ is a weak
equivalence of $\hat\phi$-squares.

(ii)\qua We show now that the composite
$\hat\Theta:=\Theta\bar\Theta\co\hat\BD\to\BD$ is a weak
equivalence in the category of commutative squares of CDGA. We
know already that $\hat\Theta$ is a quasi-isomorphism, since both
$\Theta$ and $\bar\Theta$ are. Recalling the form of $\Theta$ from
the proof of Lemma \ref{lemma-thetatheta''} and of $\bar\Theta$
from the proof of (i), we see that
$$\hat\Theta=
\left(\begin{array}{cc}%
\bar\alpha&\beta\\%
\bar\alpha\oplus s\gamma&\beta\oplus s\gamma%
\end{array}\right)
\left(\begin{array}{cc}%
\rho&\id\\%
\bar\lambda&\bar\mu%
\end{array}\right)=
\left(\begin{array}{cc}%
\alpha&\beta\\%
\lambda&\mu%
\end{array}\right),
$$
where $\alpha$, $\beta$ are CDGA morphisms and $\lambda$ (resp.\
$\mu$) is some $\hat R$-DGmodule (resp.\ $\hat Q$-DGmodule)
morphism. 

By the hypotheses of Theorem \ref{thm-stableCDGA}, we
have that $(R\oplus_{\phishrk} sD)^{>n}=0$. Since $n\geq2m+3$, this
implies that $(R\oplus_{\phishrk}sD)^{\geq2(n-m-1)}=0$. Since
$X^{<n-m-1}=0$, Lemma \ref{lemma-DGmodCDGAmap} implies that
$\lambda$ is a CDGA morphism. 

Suppose that $H^1(f)$ is injective and $n\geq 2m+3$. A similar argument shows that $\mu$
is a CDGA morphism, which implies that $\hat\Theta$ is a
weak equivalence of squares of CDGA.

Suppose instead that $n\geq 2m+4$. Since $Q$ is connected and 
$H^{\geq n}(Q\oplus sD)=0$
there exists an acyclic ideal $L\subset Q\oplus sD$ such that
$\left((Q\oplus sD)/L\right)^{\geq n}=0$. Replace $Q\oplus sD$ in diagram $\BD$
by $(Q\oplus sD)/L$ to get a quasi-)isomorphic CDGA diagram $\tilde\BD$.
Since $Y^{<n-m-2}=0$ and $2(n-m-2)\geq n$, 
Lemma \ref{lemma-DGmodCDGAmap} implies that
the composite 
$\hat Q\otimes\wedge X\otimes\wedge Y\stackrel{\mu}{\to}Q\oplus sD
\quism(Q\oplus sD)/L$
 is a CDGA quasi-isomorphism. Therefore $\hat\BD\simeq\tilde\BD\simeq\BD$ as
 CDGA squares.
\end{proof}

Collecting  these lemmas we conclude the proof of the first part of the theorem:
\begin{proof}[Proof of Theorem \ref{thm-stableCDGA}]
Lemmas \ref{lemma-thetahat''} and \ref{lemma-bartheta} imply
that the diagrams $\BD$ and $\BD'$ are weakly equivalent CDGA
commutative squares.

We prove now the second part of the theorem. 
Suppose given a CDGA model $\phi_0\co R_0\to Q_0$ of $f^*$. Our
goal is to build a model $\phi\co R\to Q$ and a \topdegree map
$\phishrk\co D\to R$ fulfilling hypotheses (i-)-(iii) of Theorem
\ref{thm-stableCDGA}. By replacing
$\phi_0$ by a minimal Sullivan model we can suppose that both
$R_0$ and $Q_0$ are connected. Since $H^{>n-1}(Q_0)=0$ and
$H^{>n}(R_0)$ we can build
another CDGA model of $f^*$ of the form
$\phi_1\co R\to Q_1$
with $R^{>n}=0$, and such that  $R$ and $Q_1$ are still connected.
We can factor $\phi_1$ into a minimal relative Sullivan algebra
$\phi_2$ followed by a weak equivalence. This gives another CDGA
model of $f^*$ of the form
$\xymatrix{
\phi_2\co R\quad\ar@{>->}[r]&Q_2:=R\otimes\wedge V}
$
and $V=V^{\geq1}$ because $H^1(f;\BQ)$ is injective.

Let $D_2$ be a minimal semifree model of the $Q_2$-DGmodule
$s^{-n}\#Q_2$. Since $H^{<n-m}(s^{-n}\#Q_2)=0$, minimality implies
that $D_2^{<n-m}=0$. Since $\phi_2$ is a relative Sullivan
algebra, every semifree $Q_2$-DGmodule is also a semifree
$R$-DGmodule. Therefore $D_2$ is also a cofibrant $R$-DGmodule and
Proposition \ref{prop-existsshriek} implies that there exists a
\topdegree map of $R$-DGmodule
$\phishrk_2\co D_2\to R$.
Since $H^{>n}(s^{-n}\#Q_2)=0$ 
we can replace $D_2$ by a weakly equivalent
$Q_2$-DGmodule, $D$, such that $D^{<n-m}=0$, $D^{>n+1}=0$, and
$D^{\leq n}=D_2^{\leq n}$.
Since $R^{>n}=0$ the map $\phishrk_2$ induces a
\topdegree map
$\phishrk\co D\to R$.
Since $H^{>m}(Q_2)=0$ and $Q_2$ is connected there exists a surjective quasi-isomorphism
of CDGA
$\alpha_2\co Q_2\quism Q$
such that $Q^{>m+2}=0$ and $\ker(\alpha_2)\subset Q^{>m+1}$. For
degree reasons $(\ker\alpha_2).D=0$, therefore the $Q_2$-DGmodule
$D$ inherits a $Q$-DGmodule structure. Set $\phi=\alpha_2\phi_2$.

In summary we have built from $\phi_0$ another CDGA model $\phi$ of $f^*$
and a \topdegree map of $R$-DGmodule
$\phishrk$ satisfying hypotheses (i)--(iii).
\end{proof}
\section{CDGA models of the complement in a Poincar\'e embedding under the unknotting condition}
\label{section-wkstableCDGA}
In this section we give a proof of Theorem \ref{thm-wkstCDGA}
which gives a CDGA model of the complement in a Poincar\'e embedding under the unknotting condition.
We also build a model of a diagram which is almost the Poincar\'e 
embedding \refequ{diag-mainsquare} under a slightly stronger unknotting condition
(Theorem \ref{thm-wkstCDGAsquare}).

The proof of Theorem \ref{thm-wkstCDGA} follows the line of the
proof of Theorem \ref{thm-stableCDGA}. In particular we will reuse
many of the lemmas of the previous section. First it is easy to
check that if we replace the hypotheses of Theorem
\ref{thm-stableCDGA} by those of Theorem \ref{thm-wkstCDGA} in
Lemmas \ref{lemma-CDGAphihat} and \ref{lemma-Rbar} then the
conclusions of these lemmas still hold without any change in their
proofs. Since $D$ is supposed to be only an $R$-DGmodule model of $s^{-n}\#Q$ we
replace Lemma \ref{lemma-stableDhat} by the following
\begin{lemma}\label{lemma-wkstDhat}
With the hypotheses of Theorem \ref{thm-wkstCDGA} and with the
notation of Lemmas \ref{lemma-CDGAphihat} and \ref{lemma-Rbar},
there exists a cofibrant $\hat R$-DGmodule, $\hat D$, and weak
equivalences of $\hat R$-DGmodules,
$$
\xymatrix{D&\ar[l]^{\simeq}_{\gamma}\hat
D\ar[r]^-{\gamma''}_-{\simeq}& \Apl(W,C),}
$$
making the following diagram of isomorphisms commute
$$\xymatrix{%
H^n(D)\ar[d]_{H^n(\phishrk)}^\cong &\ar[l]_{H^n(\gamma)}^\cong
H^n(\hat
D)\ar[r]^{H^n(\gamma'')}_\cong &H^n(W,C) \ar[d]^{H^n(\iota)}_\cong \\
H^n(R)&\ar[l]^{H^n(\bar\alpha)}_\cong H^n(\bar
R)\ar[r]_{H^n(\bar\alpha')}^\cong & H^n(W).}
$$
\end{lemma}
\begin{proof}
It is a special case of Lemma \ref{lemma-DGmodDhat}.
\end{proof}
It can be readily checked that Lemma
\ref{lemma-CDGAphitildeshriek} still holds when we replace the
hypotheses of Theorem \ref{thm-stableCDGA} by those of Theorem
\ref{thm-wkstCDGA}, and the only change in the proof of this lemma
is a replacement of the reference to Lemma \ref{lemma-stableDhat}
to a reference to Lemma \ref{lemma-wkstDhat}.
Moreover by Lemma \ref{lemma-rigidify} we can replace 
$\tilde\phishrk$ by $\bar\phishrk$
making the diagram of Lemma
\ref{lemma-CDGAphitildeshriek} strictly commute.
We are now ready for the following
\begin{proof}[Proof of Theorem \ref{thm-wkstCDGA}]
By the same argument as for Theorem \ref{thm-stableCDGA} and using 
Lemma \ref{lemma-CDGAphihat}, \ref{lemma-Rbar},
\ref{lemma-wkstDhat}, and \ref{lemma-CDGAphitildeshriek}, we get that
 $R\hookrightarrow R\oplus_{\phishrk}sD$ is an $\hat R$-DGmodule
model of $l^*\co\Apl(W)\hookrightarrow\Apl(C)$.

Since $H_*(f)$ is $r$-connected and by Lefschetz duality we have
 $H^{\geq n-r}(R\oplus_{\phishrk}sD)=H^{\geq n-r}(C;\BQ)=H_{\leq r}(W,P;\BQ)=0$.
 Using the connectivity of $R$ it is easy to build an acyclic subDGmodule
 $L\subset R\oplus_{\phishrk}sD$ concentrated in degrees $\geq n-r-1$ and
 killing $(R\oplus_{\phishrk}sD)^{\geq n-r}$.

Consider such an acyclic subDGmodule $L$. Set $k=n-m-1$ and $l=n-2m+r-1$. 
By  Lemma \ref{lemma-CDGAtruncMC}
there is a semi-trivial CDGA structure on
$(R\oplus_{\phishrk}sD)/L$ and the obvious map $R\to
(R\oplus_{\phishrk}sD)/L$ is a CDGA morphism.
Since $L$ is acyclic the map $R\to (R\oplus_{\phishrk}sD)/L$
is also an $\hat R$-DGmodule model $l^*$.

Let $\xymatrix@1{ \hat R\quad\ar@{>->}[r]^-u&\hat R\otimes\wedge
X\ar[r]^{\lambda'}_{\simeq}&\Apl(C)}$
 be a CDGA factorization of $l^*\alpha'$
through a relative minimal Sullivan algebra $\hat R\otimes\wedge
X$. By the same argument as in the proof of Lemma
\ref{lemma-thetahat''} we find that $X^{<n-m-1}=0$ because $l^*$
is $(n-m-1)$-connected.

Since $u$ is a model of $\hat R$-DGmodule of $l^*$, the same
argument as in the beginning of the proof of Lemma
\ref{lemma-bartheta} gives a commutative diagram of $\hat
R$-DGmodules
$$
\xymatrix{ \hat R_{{}_{{}_{{}_{}}}}\ar[r]^{\rho}_{\simeq}\ar@{>->}[d]_u&\bar
R_{{}_{{}_{{}_{}}}}\ar@{^(->}[d]\ar[r]^{\bar\alpha}_{\simeq}&
R_{{}_{{}_{{}_{}}}}\ar[d]\\
\hat R\otimes\wedge X\ar[r]^{\bar\lambda}_\simeq
\ar@(d,d)[rr]^{\lambda}_{\simeq}&\bar R\oplus_{\bar\phishrk}s\hat
D\ar[r]^-{\pi(\bar\alpha\oplus
s\gamma)}_-\simeq&(R\oplus_{\phishrk}sD)/L,}
$$
 and the composite $\lambda=\pi(\bar\alpha\oplus
s\gamma)\bar\lambda$ is a quasi-isomorphism.

Since $X^{<n-m-1}=0$ and $\left((R\oplus_{\phishrk}sD)/L\right)^{\geq n-r}=0$, 
the condition $r\geq2m-n+2$ and Lemma
\ref{lemma-DGmodCDGAmap} imply that $\lambda$ is a CDGA morphism.
 Also $\bar\alpha\rho=\alpha$ is a CDGA morphism. Thus $u$ is
a CDGA model of $R\to(R\oplus_{\phishrk}sD)/L$. By
construction $u$ is also a model of $l^*$ and the first part of the theorem is
proved.

The second part of the theorem is proved in a similar way to Theorem \ref{thm-stableCDGA}.
\end{proof}

\begin{proof}[Proof of Corollary \ref{corol-wkstCDGA}]
Since $P$ is simply connected and the codimension is at least $3$,
$\del T$ is simply connected, and since $W$ is also simply-connected, the same is true for $C$ by 
 Van Kampen theorem. The corollary follows then from the above theorem
 and from the fact that a CDGA model of a  simply connected
spaces of finite type determines its rational homotopy type.
\end{proof}

In the rest of the section we address the problem of describing a
CDGA model of Diagram \refequ{diag-mainsquare} under some
unknotting condition.
We wish that we could have determined  the
rational homotopy type of the entire square
\refequ{diag-mainsquare} from the rational homotopy class of $f$,
but we are only able to determine a slightly less complete  square
that we describe now. 

Assume that $\del T$ is simply-connected in which case
by Poincar\'e duality in dimension $n-1$ and by
\cite[Proposition 4.1]{Wall-finiteness} we can consider the
space $\del\punct T$ obtained by removing the unique top $(n-1)$-cell in
 a minimal CW-decomposition of $\del T$. We have then the following
commutative square of topological spaces
\begin{equation}\label{diag-mainsquarepunct}%
 \xymatrix{ \del\punct T\ar[r]^{\punct i}\ar[d]_{\punct k}&
P\ar[d]^f\\
C\ar[r]_l&W }
\end{equation}
where  $\punct i$ and $\punct k$ are the
restrictions of $i$ and $k$ to $\del \punct T$. Our next theorem
is a description of a CDGA model of \refequ{diag-mainsquarepunct}
under a stronger unknotting condition  and two
extra assumptions which are not too restrictive as we explain in
Remark \ref{rmk-extraassumption}.

To state the theorem it is convenient to introduce the following terminology:
 if $X$ is an $A$-DGmodule and $l$ is an integer then a \emph{truncation
$A$-subDGmodule of $X$ above degree $l$}  is a subDGmodule $L$ such that $L^{\leq l-1}=0$,
$X^{>l}\subset L$ and the projection $\pi\co X\to X/L$ induces an isomorphism
in homology in degrees $\leq l$. Of course $(X/L)^{>l}=0$. It is easy to check that such a truncation
subDGmodule exists when $A$ is connected.
\begin{thm}\label{thm-wkstCDGAsquare}
Consider the diagram \refequ{diag-mainsquarepunct} induced
by a Poincar\'e embedding \refequ{diag-mainsquare} with
$P$ and $W$ connected and $\del T$ simply-connected. Let $r$ be a positive integer
such that $\tilde H_{\leq r-1}(P;\BQ)=\tilde H_{\leq r}(W;\BQ)=0$
and $r\geq 2m-n+2$.
Let $\phi\co R\to Q$ be a CDGA model of
$f^*\co\Apl(W)\to\Apl(P)$ such that $R$ is connected. Let D be
an $R$-DGmodule weakly equivalent to $s^{-n}\#Q$ and such that
$D^{<n-m}=0$. Suppose given a \topdegree map of $R$-DGmodules
$\phishrk\co D\to R$. 

Suppose moreover that $n\geq m+r+2$ and that $Q$
is $(r-1)$-connected, that is $Q^{\leq r-1}=Q^0=\BQ$.

Let $I$ be a truncation $R$-subDGmodule of $R$ above degree $n-r-1$,
let $J$ be a truncation $Q$-subDGmodule of $Q$ above degree $m$, and 
let $K$ be a truncation $Q$-subDGmodule of $D$ above degree $n-r$.

Then the following two commutative squares are weakly equivalent
in CDGA
$$\BDt:=
 \vcenter{\xymatrix{ R\ar[r]^\phi\ar[d]& Q\ar[d]\\
(R\oplus_{\phishrk}sD)/(I\oplus
sK)\ar[r]^{\overline{\phi\oplus\id}}&
(Q\oplus_{\phi\phishrk} sD)/(J \oplus sK)}}
$$
and
$$
\BDt':=%
\vcenter{\xymatrix{%
 \Apl(W)\ar[r]^{f^*}\ar[d]_{l^*}&
\Apl(P)\ar[d]^{\punct i^*}\\
\Apl(C)\ar[r]^{\punct k^*}& \Apl(\del\punct T). }}
$$
where, in Diagram $\BDt$, the vertical maps are the composition of
the inclusion with the projection, the bottom map is the one
induced by $\phi\oplus\id_{sD}$, and the CDGA structure on the
truncated mapping cones are the semi-trivial ones.
\end{thm}
\begin{rmk}\label{rmk-extraassumption}
The connectivity hypothesis on $P$ and $W$ are equivalent 
to $H_*(f;\BQ)$ is $r$-connected and $\tilde H_{\leq r-1}(P;\BQ)=0$
which is clearly a stronger condition than the unknotting condition \refequ{equ-unknotrht}
because of the high connectivity hypothesis on $P$.
The first extra assumption in the theorem, $n\geq
m+r+2$, is satisfied under the unknotting condition $r\geq 2m-n+2$ as
soon as $m\geq 2r$. On the other hand, if $m<2r$ then by a
rational version of the suspension  Freudenthal theorem, $P$ has
the rational homotopy type of a wedge of spheres of dimensions
between $r$ and $2r-1$. Hence this first extra assumption is a
consequence of the unknotting condition when $P$ is not rationally
equivalent to a wedge of spheres. For the second extra assumption
(the $(r-1)$-connectivity of $Q$), since $\tilde H^{\leq r}(P)=0$, one
can always construct an $r$-connected CDGA model $Q$ of $P$, by
taking for example a minimal Sullivan model of any given model of
$P$. Therefore there is no real loss of generality in making this second
assumption.
\end{rmk}
\begin{rmk} It is very likely that under the only unknotting condition \refequ{equ-unknotrht} 
one can determine a CDGA model of the
complete Poincar\'e embedding \refequ{diag-mainsquare} but we were
unable to prove this.
\end{rmk}

The rest of the section is devoted to the proof of Theorem
\ref{thm-wkstCDGAsquare} which is a refinement of the proof of
Theorem \ref{thm-stableCDGA}.

Lemmas \ref{lemma-CDGAphihat}--\ref{lemma-CDGAphitildeshriek} 
hold with the
hypotheses of Theorem \ref{thm-stableCDGA} replaced by those of
Theorem \ref{thm-wkstCDGAsquare}.
We need the following three lemmas in replacement of Lemma
\ref{lemma-chi}:
\begin{lemma}
\label{lemma-chiwkst} With the hypotheses of Theorem
\ref{thm-wkstCDGAsquare} and with the notation of Lemmas
\ref{lemma-CDGAphihat}--\ref{lemma-CDGAphitildeshriek}, the
composite $\phi\phishrk$ induces a $Q$-DGmodule map
$
\overline{\phi\phishrk}\co D/K\to Q/J$. There exists a $\hat
Q$-DGmodule morphism $\chi\co \hat D\to\hat Q$
making the following diagram commute in $\hat Q$-DGMod
$$
\xymatrix{
D/K\ar[d]_{\overline{\phi\phishrk}}&\ar[l]_{\pi_0}D&\ar[l]^{\simeq}_{\gamma}\hat
D\ar[r]^-{\gamma'}_-{\simeq}\ar[d]^\chi&\Apl(P,\del T)\ar[d]^{\iota'}\\
Q/J&\ar[l]_{\pi'_0}Q&\ar[l]^{\simeq}_{\beta}\hat
Q\ar[r]^{\beta'}_{\simeq}&\Apl(P)}
$$
where $\pi_0$ and $\pi'_0$ are the canonical projections. Moreover
$\chi\simeq_{\hat R}\bar\phi\tilde\phishrk$.
\end{lemma}
\begin{proof}
 Since $n\geq m+r+2$ we have that $\phi\phishrk(K)\subset J$,
hence there is an induced map $\overline{\phi\phishrk}$ between the
quotients. Since $Q$ is $(r-1)$-connected, $(D/K)^{<n-m}=0$,
$(Q/J)^{>m}=0$, and $(n-m)+r>m$, we have
that the $\Bk$-DGmodule map $\overline{\phi\phishrk}$ is a map of
$Q$-DGmodule.

Since $r\geq2m-n+2$ and $D^{<n-m}=0$, we have $H^{\leq m-r+1}(\hat D)=0$.
Also $\tilde H^{\leq r-1}(\hat Q)=H^{>m}(\hat Q)=0$. By Lemma
\ref{lemma-stablehmtpyclasses} we have an isomorphism
$
 H^*\co[\hat D,\hat Q]_{\hat
Q}\cong\homset^0_{\Bk}(H^*(\hat D),H^*(\hat Q))
$.
Therefore there exists a map
$\chi\co\hat D\to\hat Q$ of $\hat Q$-DGmodules, unique up to
homotopy, such that $H^*(\chi)=H^*(\bar\phi\tilde\phishrk)$ where
$\bar\phi$ and $\tilde\phishrk$ were defined in Lemma
\ref{lemma-Rbar} and Lemma \ref{lemma-CDGAphitildeshriek}.

Since $\chi$ induces in cohomology the same map as
$\bar\phi\tilde\phishrk$, Lemma \ref{lemma-CDGAphitildeshriek},
Lemma \ref{lemma-CDGAphihat} and Lemma \ref{lemma-Rbar} imply that
the map $\chi$ makes the diagram of the statement of  Lemma
\ref{lemma-chiwkst} commute \emph{in cohomology}. Another
application of Lemma \ref{lemma-stablehmtpyclasses} implies that
this diagram commutes up to a homotopy of $\hat Q$-DGmodules.
Since $(\beta\pi'_0,\beta')$ is surjective we can suppose by Lemma
\ref{lemma-rigidify} that $\chi$ makes the diagram exactly
commute.

Finally we have also $\chi\simeq_{\hat R}\bar\phi\tilde\phishrk$, again by Lemma
\ref{lemma-stablehmtpyclasses}.
\end{proof}

\begin{lemma}\label{lemma-connectivityHconebis}
 With the hypotheses of Theorem
\ref{thm-wkstCDGAsquare}, the canonical projection
$$\pi\co R\oplus_{\phishrk} sD\to (R\oplus_{\phishrk} sD)/(I\oplus sK)$$
is a quasi-isomorphism and the
canonical projection
$$\pi'\co Q\oplus_{\phi\phishrk} sD\to (Q\oplus_{\phi\phishrk} sD)/(J\oplus sK)$$
induces an isomorphism in cohomology in all degrees except in
degree $n-1$ where $H^{n-1}(Q\oplus_{\phi\phishrk} sD)\cong\BQ$ and
$H^{n-1}((Q\oplus_{\phi\phishrk} sD)/(J\oplus sK))=0$.
\end{lemma}
\begin{proof}
$L=I\oplus sK$ is a truncation $R$-subDGmodule of $R\oplus_{\phishrk} sD$
above degree $n-r$. Using the fact that $H^n(\phishrk)$ is an isomorphism
and that $H^i(R)=H^i(sD)=0$ for $n-r\leq i\not =n$,
it comes that $L$ is acyclic, hence $\pi$ is a quasi-isomorphism.

The proof for $\pi'$ is similar after computing that
$H^{\geq n-r}(Q\oplus_{\phi\phishrk} sD)\cong s^{-(n-1)}\BQ$ and using
the  assumption $n\geq m+r+2$ to check that $J\oplus sK$ is a differential submodule
of $Q\oplus_{\phishrk}sD$.
\end{proof}
\begin{lemma}\label{lemma-modeldTcheck}
 With the hypotheses of Theorem \ref{thm-wkstCDGAsquare} and
 with the notation of Lemmas
\ref{lemma-CDGAphihat}--\ref{lemma-CDGAphitildeshriek} and Lemma
\ref{lemma-chiwkst}, there exists a cofibration of $\hat
Q$-DGmodules
$$\xymatrix{
w\co\hat Q\oplus_\chi s\hat D\,\,\ar@{>->}[r] &(\hat
Q\oplus_\chi s\hat D)\oplus \hat Q\otimes V}
$$
and acyclic fibrations of $\hat Q$-DGmodules $\epsilon$ and
$\epsilon'$ making the following diagram commute
$$\xymatrix{%
Q\oplus_{\phi\phishrk}sD\ar@{->>}[d]_{\pi '}&%
\ar@{->>}[l]_{\alpha\oplus s\gamma}^\simeq {\hat Q\oplus_\chi s\hat
D}_{{}_{{}_{{}_{}}}}\ar@{>->}[d]_w\ar@{->>}[r]^{i^*\alpha'\oplus 0}_\simeq&%
\Apl(\del T)\ar@{->>}[d]^{\check j}\\
(Q\oplus_{\phi\phishrk}sD)/(J\oplus sK)&%
\ar@{->>}[l]_-{\epsilon}^-\simeq (\hat Q\oplus_\chi s\hat
D)\oplus \hat Q\otimes V\ar@{->>}[r]^-{\epsilon'}_-\simeq&%
\Apl(\del \check T). }
$$
\end{lemma}
\begin{proof}
The composite of $\hat Q$-DGmodules
$$
\xymatrix{%
\hat Q\oplus_\chi s\hat D%
\ar@{->>}[r]^{i^*\alpha'\oplus 0}_\simeq&%
\Apl(\del T)\ar@{->>}[r]^{\check j}&%
\Apl(\del \check T)}
$$
can de factored into a \emph{minimal} semi-free extension $w$
followed by a quasi-isomorphism $\epsilon'$. Moreover, since
$\check j(i^*\alpha'\oplus 0)$ is a surjection, so is $\epsilon'$.

Define $\epsilon$ as the extension of $\pi'(\alpha\oplus s\gamma)$ such
that $\epsilon(\hat Q\otimes V)=0$, which is a $\hat Q$-module morphism.
It is clear that $\check j$ is $(n-2)$-connected and by minimality
$V^{<n-2}=0$. Since $r\geq 1$, we have
$((Q\oplus sD)/(J\oplus sK))^{\geq n-1}=0$.
For degree reasons $\epsilon$ is a DGmodule map.

It remains to prove that $\epsilon$ is a quasi-isomorphism. This
is an easy consequence of the fact that $H^{<n-1}(\pi')$ is an
isomorphism and $H^{\geq n-1}((Q\oplus_{\phi\phishrk}sD)/(J\oplus
sK))=0=H^{\geq n-1}(\partial \check T)$.
\end{proof}

Lemma \ref{lemma-factorchi} holds with the hypotheses of Theorem
\ref{thm-wkstCDGAsquare} replacing those of Theorem
\ref{thm-stableCDGA}, without any change in the proof.

To finish the proof of Theorem \ref{thm-wkstCDGAsquare}, we  adapt
the four Lemmas \ref{lemma-phihatsq}--\ref{lemma-bartheta} to the
setting of this section. Recall Diagrams $\BDt$ and $\BDt'$ from
the statement of Theorem \ref{thm-wkstCDGAsquare}.

\begin{lemma}\label{lemma-wkphihatsq}
With the hypotheses of Theorem \ref{thm-wkstCDGAsquare} and  with
the notation of Lemmas
\ref{lemma-CDGAphihat}--\ref{lemma-factorchi} and
\ref{lemma-chiwkst}--\ref{lemma-modeldTcheck}, Diagrams $\BDt$ and
$\BDt'$ are both commutative squares in CDGA and
$\hat\phi$-squares,  and the following diagram is a
$\hat\phi$-square:
$$
\bar \BDt:=%
\vcenter{\xymatrix@1{%
\bar R\vrule width 0pt depth 6pt\ar[r]^{\bar\phi}\ar@{^(->}[d]&%
\hat Q\vrule width 0pt depth 6pt\ar@{^(->}[d]\\%
\bar R\oplus_{\bar\phishrk}s\hat D\ar[r]_-{w(\bar\phi\oplus\id)}&%
(\hat Q\oplus_\chi s\hat D)\oplus \hat Q\otimes V%
}}%
$$
\end{lemma}
\begin{proof}
We show first that the $\BDt$ is a diagram of CDGA. We have
already shown in the proof of Theorem \ref{thm-wkstCDGA} that
$R\to (R\oplus_{\phishrk}sD)/(I\oplus sK)$ is a CDGA map. The
morphism $\phi\phishrk$ is not a $Q$-DGmodule morphism but, for
degree reasons, the composite $\pi'_0\phi\phishrk\co D\to Q/J$
is. Therefore the truncated mapping cone $(Q\oplus_{\phi\phishrk}
sD)/(J\oplus sK))$ has a natural structure of $Q$-DGmodule. Again
by Lemma \ref{lemma-CDGAtruncMC}, this endows this mapping cone
with a semi trivial CDGA structure and the map $Q\to
(Q\oplus_{\phi\phishrk} sD)/(J\oplus sK)$ is a CDGA map. Moreover,
using the fact that $n\geq m+r+2$ we get that $\phi(I)\subset J$,
therefore $\phi\oplus \id\co R\oplus_{\phishrk}sD\to Q
\oplus_{\phi\phishrk}sD$ induces a map, $\overline{\phi\oplus \id}$,
between the quotients. It is straightforward to check that it is a
CDGA map.

This proves that $\BDt$ is a CDGA square and also a
$\hat\phi$-square where the $\hat R$- and $\hat Q$-DGmodule
structures are induced by the maps $\alpha$ and $\beta$. It is
immediate that $\BDt'$ is a CDGA-square and it is also a
$\hat\phi$-square where the $\hat R$- and $\hat Q$-DGmodule
structures are induced by the maps $\alpha'$ and $\beta'$.

It is immediate that $\bar\BDt$ is  a $\hat\phi$-square.
\end{proof}
\begin{lemma}\label{lemma-wkthetatheta''}
With  the hypotheses of Theorem \ref{thm-wkstCDGAsquare} and  with
the notation of Lemmas
\ref{lemma-CDGAphihat}--\ref{lemma-factorchi} and
\ref{lemma-chiwkst}--\ref{lemma-wkphihatsq}, there exist acyclic
fibrations of $\hat\phi$-squares
$$
\xymatrix{\BDt&%
\ar@{->>}[l]_{\Theta}^{\simeq}\bar\BDt\ar@{->>}[r]^{\Theta'}_\simeq&%
\BDt'.}
$$
\end{lemma}
\begin{proof}
Set
$\Theta=\left(\begin{array}{cc}\bar\alpha&\beta\\\pi(\alpha\oplus
s\gamma)&\epsilon\end{array}\right)$ and
$\Theta'=\left(\begin{array}{cc}\bar\alpha'&\beta'\\\l^*\bar\alpha'\oplus
0&\epsilon'\end{array}\right)$ where $\epsilon$ and $\epsilon'$
were defined in Lemma \ref{lemma-modeldTcheck}. An argument
analogous to that of Lemma \ref{lemma-thetatheta''} together with
the results of Lemmas \ref{lemma-connectivityHconebis} and
\ref{lemma-chiwkst} finishes the proof.
\end{proof}

\begin{lemma}\label{lemma-wkthetahat''}
With the hypotheses  of Theorem \ref{thm-wkstCDGAsquare} and  with
the notation of Lemmas
\ref{lemma-CDGAphihat}--\ref{lemma-factorchi} and
\ref{lemma-chiwkst}--\ref{lemma-wkthetatheta''}, there exists a commutative square in CDGA,
$$\hat\BDt:=%
\vcenter{\xymatrix@1{%
\vrule width0pt depth8pt\hat R_{{}_{{}_{{}_{}}}}\quad\ar@{>->}[r]^{\hat\phi}\ar@{>->}[d]_u&%
\vrule width0pt depth8pt\hat Q_{{}_{{}_{{}_{}}}}\ar@{>->}[d]^v\\
\hat R\otimes\wedge X\,\,\ar@{>->}[r]_-{\hat\psi}&%
\hat Q\otimes \wedge X\otimes\wedge Z}}$$
where $\hat\phi$, $\hat\psi$, $u$, and $v$  are cofibrations, and
there exists a weak equivalence both of CDGA-squares and of
$\hat\phi$-squares
$\hat\Theta'\co\hat\BDt\quism \BDt'$.
Moreover $X$ and $Z$ can be chosen such that such that
$X^{<n-m-1}=Z^{<n-m-1}=0$.
\end{lemma}
\begin{proof}
The proof is completely similar to that of Lemma
\ref{lemma-thetahat''}, replacing $\Apl(\del T)$ by $\Apl(\del
\check T)$, which changes nothing to the $(n-m-1)$-connectivity of
the maps and noticing that since $r$ is positive, $H^1(f)$ is injective.
\end{proof}

\begin{lemma}\label{lemma-wkbartheta}
With the hypotheses  of Theorem \ref{thm-wkstCDGAsquare} and  with
the notation of Lemmas
\ref{lemma-CDGAphihat}--\ref{lemma-factorchi} and
\ref{lemma-chiwkst}--\ref{lemma-wkthetahat''}, there exists a
quasi-isomorphism of commutative squares in CDGA
$\hat\Theta\co\hat\BDt\quism\BDt$.
\end{lemma}
\begin{proof}
 By an completely  analogous
argument to that of the first part of the proof of Lemma
\ref{lemma-bartheta}, we get a lifting of $\hat\phi$-squares
$\bar\Theta\co\hat\BDt\quism\bar\BDt$.
It remains then to prove that the composite
$\hat\Theta:=\bar\Theta\hat\Theta'$ is a morphism of squares of
CDGA. This is proved by Lemma \ref{lemma-DGmodCDGAmap} using the
facts that $$((R\oplus_{\phishrk}sD)/(I\oplus sK))^{\geq n-r}=
((Q\oplus_{\phi\phishrk}sD)/(J\oplus sK))^{\geq n-r}=0,
$$
that $X^{<n-m-1}=Z^{<n-m-1}=0$, and that $2(n-m-1)\geq n-r$ by
\refequ{equ-unknotrht}.
\end{proof}

\begin{proof}[Proof of Theorem \ref{thm-wkstCDGA}]
The fact that $\BDt$ is a well defined CDGA square was proved in
the first part of Lemma \ref{lemma-wkphihatsq}. Then by Lemmas
\ref{lemma-wkthetahat''} and \ref{lemma-wkbartheta}, the diagrams
$\BDt$ and $\BDt'$ are weakly equivalent in CDGA.
\end{proof}
\section{Examples of rationally knotted embeddings}
\label{section-examples} The aim of this section is to show by
some examples that the unknotting condition
\refequ{equ-unknotrht}
in Theorem \ref{thm-wkstCDGA} and in the second part of 
Corollary \ref{corol-HWmodstruct} is unavoidable and sharp. Recall
that this condition is $r\geq 2m-n+2$ where
\begin{itemize}
\item $m$ is the dimension of the embedded
polyhedron $P$,
\item $n-m$ is the codimension of the embedding,
\item $r$ is the connectivity of the embedding.
\end{itemize}

We will build two families of examples where the unknotting
condition \refequ{equ-unknotrht} is missed by a little and such
that the thesis of Theorem \ref{thm-wkstCDGA} does not hold. The
unknotting condition can be reformulated as $r+(n-m)\geq m+2$ which
 can be roughly expressed as $$
\mathrm{connectivity}+\mathrm{codimension}\,\geq\,\mathrm{dimension}+2.
$$
In the first examples that we will build (Proposition
\ref{prop-exsharp1}), the connectivity $r$ is big but the
codimension $n-m$ is not high enough, and in the second family of
examples (Proposition \ref{prop-exsharp2}) the codimension will be
big but the connectivity small. Both of these families of examples
are fairly explicit and  are described in the proof of these
propositions.

\begin{prop}
\label{prop-exsharp1} Let $p$ be a positive even integer and let
$n\geq 3p+2$. Set $m=n-p-1$ and $r=2m-n+1$.
Then there exist two  $m$-dimensional polyhedra,
$P_0$ and $P_1$, having both the rational homotopy type of the wedge of spheres
$S^{n-2p-1}\vee S^{n-p-1}$, and two
nullhomotopic $r$-connected embeddings $f_0\co P_0\hookrightarrow S^n$ and
$f_1\co P_1\hookrightarrow S^n$, such that the rational
cohomology algebras of the complement of these embedded polyhedra
are not isomorphic:
$$H^*(S^n\takeaway f_0(P_0);\BQ)\not\cong
H^*(S^n\takeaway f_1(P_1);\BQ).
$$
\end{prop}
\begin{proof}
Set $X_0=S^p\vee S^{2p}$. There exists an obvious PL-embedding
$X_0\subset S^n$. Define $P_0$ as the closure of the complement of
some regular neighborhood of $X_0$ in $S^n$. By Lefschetz duality we have
$\tilde H_*(P_0;\BZ)=\BZ.x_{n-2p-1}\oplus \BZ.y_{n-p-1}$ and by
\cite[Proposition 4.1]{Wall-finiteness}
$P_0$ has the homotopy type of a  two-cell CW-complex
$
P_0 \simeq S^{n-2p-1}\cup e^{n-p-1}$. Since $n\geq 3p+2$, we have
that $\pi_{n-p-2}(S^{n-2p-1})\otimes\BQ=0$ and $ P_0\simeq_\BQ S^{n-2p-1}\vee S^{n-p-1}$.
Therefore the rational cohomology algebra $
H^*(S^n\takeaway P_0,\BQ)\cong H^*(X_0;\BQ)$
has a trivial multiplication.

On the other hand  consider a $(p-1)$-connected and
$2p$-dimensional polyhedron $X_1$ having the homotopy type of the
CW-complex $S^p\cup_{2[\iota,\iota]}e^{2p}$
where $\iota\in\pi_p(S^p)$ represents the identity map and
$[\iota,\iota]$ is the Whitehead bracket. Then
$H^*(X_1;\BQ)\cong \BQ[x]/(x^3)$ with $\deg(x)=p$. By the
embedding  theorem of Wall \cite{Wall-thick}, after replacing
$X_1$ by some polyhedron of the same homotopy type, there exists
an embedding $X_1\subset S^n$. Define $P_1$ as the closure of the
complement of a regular neighborhood of $X_1$ in $S^n$. By the
same argument as for $P_0$ we see that $P_1$ has the rational
homotopy type of the same wedge of spheres.
But here the multiplication on the cohomology algebra
$H^*(S^n\takeaway P_1;\BQ)\cong H^*(X_1;\BQ)$
is \emph{not} trivial.

Finally it is immediate that both embeddings $P_0\subset S^n$ and
$P_1\subset S^n$ are nullhomotopic and $r$-connected.
\end{proof}

The previous proposition implies that there is no way of getting a
model of the rational homotopy type of the complement
$S^n\takeaway f_i(P_i)$ from just a model of the homotopy class of
the embedding $f_i$. Notice that the equation $r=2m-n+1$ is very
close to the unknotting condition \refequ{equ-unknotrht}, showing
that this condition  is sharp in Theorem \ref{thm-wkstCDGA}
and Corollary \ref{corol-HWmodstruct}.

Note also that using Spanier-Whitehead duality or the techniques
of \cite{LSV-TAMS}, it can be shown that the two polyhedra $P_0$
and $P_1$ constructed in Proposition \ref{prop-exsharp1} can be
chosen as having the same integral homotopy type. It is even
possible that $P_0$ and $P_1$ might be chosen as being
PL-homeomorphic, but we have no proof of that fact.\medbreak
 The
examples of Proposition \ref{prop-exsharp1} show that the
unknotting condition of Theorem \ref{thm-wkstCDGA} is sharp at
least when the codimension $n-m$ is low (even if the connectivity
$r$ is high). In the rest of this section we will build a second
family of examples  for which the codimension is high but the
connectivity is low. We prove first a lemma.
\begin{lemma}
\label{lemma-embsusp} Let $i\co X\hookrightarrow S^{n-1}$ be
the inclusion of a polyhedron in a sphere and denote by
$\epsilon\co S^{n-1}\hookrightarrow S^n$ the inclusion of the
equator. Then
$$
S^n\takeaway\epsilon(i(X))\simeq \Sigma\left(S^{n-1}\takeaway
i(X)\right).
$$
\end{lemma}
\proof
Set $Y=S^{n-1}\takeaway i(X)$. It is clear that the complement of
$X$ in $S^n$ has the homotopy type of two disks $D^n$ glued along
$Y\subset S^{n-1}=\partial D^n$. Thus
$$
S^n\takeaway i(X)\simeq D^n\cup_Y D^n\simeq \Sigma Y.
\eqno{\qed}$$

\begin{prop}
\label{prop-exsharp2} For $0\leq r\leq 5$ there exists two $r$-connected
homotopic embeddings $f_k\co S^r\times S^7\hookrightarrow
S^{15}$, $k=0,1$, such that the rational cohomology algebras of
their complement are not isomorphic,
$$
H^*(S^{15}\takeaway f_0(S^r\times S^7);\BQ)\not\cong
H^*(S^{15}\takeaway f_1(S^r\times S^7);\BQ).
$$
\end{prop}
\begin{proof}
We have the standard embeddings $S^r\subset \BR^{r+1}$ and
$S^7\subset \BR^8$, as well as the ``stereographic'' embedding
$\BR^{r+9}\subset (\BR^{r+9}\cup\set{\infty})\cong S^{r+9}$.
Composing those we get an embedding
$$
i\co S^r\times S^7\hookrightarrow
\BR^{r+1}\times\BR^8=\BR^{r+9}\hookrightarrow S^{r+9}.
$$
Since $r+9\leq 14$, we have the inclusion of a subequator
$$
\epsilon\co S^{r+9}\subset S^{15}.
$$
Set $f_0=\epsilon i$. Lemma \ref{lemma-embsusp} implies that
$S^{15}\takeaway f_0(S^r\times S^7)$ has the homotopy type of a
suspension. Therefore the multiplication on $H^*(S^{15}\takeaway
f_0(S^r\times S^7);\BQ)$ is trivial.

We construct now another embedding $f_1$. Consider the Hopf
fibration
$$
S^7\to S^{15}\stackrel{\pi}{\to}S^8.$$ Consider the inclusion of
$S^r$ in $S^8$ as a subequator. Its complement $S^8\takeaway S^r$
has the homotopy type of $S^{7-r}$. Therefore the sphere $S^{15}$
is the union of two polyhedra of the homotopy type of
$\pi^{-1}(S^r)$ and $\pi^{-1}( S^{7-r})$. Since both of the
inclusions $S^r\subset S^8$ and $S^{7-r}\subset S^8$ are
nullhomotopic, the restrictions of the Hopf fibration to these
subspaces are trivial, hence $\pi^{-1}(S^r)\simeq S^r\times S^7$
and $\pi^{-1}(S^{7-r})\simeq S^{7-r}\times S^7$.
 This
defines an embedding $f_1\co S^r\times S^7\hookrightarrow
S^{15}$ whose complement has the homotopy type of $S^{7-r}\times
S^7$. Therefore the multiplication on the cohomology algebra $
H^*(S^{15}\takeaway f_1(S^r\times S^7);\BQ)$ is not trivial.

Finally it is immediate that the embeddings $f_0$ and $f_1$ are
homotopic since there are both nullhomotopic for
dimension-connectivity reasons.
\end{proof}
Taking $r=0$ in Proposition \ref{prop-exsharp2} gives an example
of two homotopic $0$-connected embeddings of $S^0\times S^7$ in $S^{15}$, of
relatively high codimension, and whose complement do not have the
same rational homotopy type. Again this shows that the unknotting
condition \refequ{equ-unknotrht} is sharp since here $r=2m-n+1$. 
Note that $r=0$ is not a positive integer and $P=S^0\times S^7$ is not connected 
as it should be in the hypotheses of Theorem
\ref{thm-wkstCDGA}. But if we take $r=1$ we get
two $1$-connected homotopic embeddings of $S^1\times S^7$ into $S^{15}$, and the
unknotting condition is only missed by $2$ in that case.

Examples analogous to those of Proposition \ref{prop-exsharp2} can
be build in other dimensions by replacing the Hopf fibration 
$S^7\to S^{15}\to S^8$ by
the Stiefel fibration
$$
S^{2k-1}\to V_2(\BR^{2k+1})\stackrel{\pi'}{\to}S^{2k}
$$
where $V_2(\BR^{2k+1})$ can be seen as the spherical tangent
bundle of $S^{2k}$. Since the Euler characteristic of an
even-dimensional sphere is not zero, it is immediate that
$V_2(\BR^{2k+1})$ has the rational homotopy type of a sphere
$S^{4k-1}$. We leave to the reader the details of the statement
and proof of a proposition analogous to \ref{prop-exsharp2} with
two embeddings of $S^r\times S^{2k-1}$ into
$V_2(\BR^{2k+1})\simeq_\BQ S^{4k-1}$ for which the rational
cohomology algebras of the complements are not isomorphic.

\Addresses\recd

\begin{thebibliography}
\bibitem{Boilley}
{\bf C. Boilley},
{\it Extensions de modules et dualit\'e de Lefschetz},
In preparation.

\bibitem{BG}
{\bf A.K. Bousfield, V.K.A.M. Gugenheim}, {\it On PL DeRham
theory and rational homotopy type}, Mem. Amer. Math. Soc. 8 (1976) 179
\MR{0425956O}



\bibitem{Bredon}
{\bf G.E.~Bredon}, {\it Topology and Geometry}, Corrected third
printing, Graduate Texts in
Mathematics 139, Springer--Verlag (1997)
\MR{1700700}


\bibitem{Deligne}
{\bf P. Deligne}, {\it Th\'eorie de Hodge, II}, Publ. Math.
I.H.E.S. {40} (1971) 5--58
  \MR{0498551}


\bibitem{DwyerSpalinski}
{\bf W.G.~Dwyer, J.~Spalinski}, {\it Homotopy theories and model
categories}, In: Handbook of Algebraic Topology, edited by
I.M.~James, Elsevier Science B.V. (1995) 73--126
  \MR{1361887}

\bibitem{FHT-gorenstein}
{\bf Y.~F\'elix, S.~Halperin, J.C.~Thomas}, {\it Gorenstein
spaces}, Advances in Math {71} (1988) 92--112
  \MR{0960364}

\bibitem{FHT-dgatop}
{\bf Y.~F\'elix, S.~Halperin, J.C.~Thomas}, {\it Differential
graded algebras in topology}, In: Handbook of Algebraic Topology,
edited by I.M.~James, Elsevier Science B.V. (1995) 829--865
  \MR{1361901}

\bibitem{FHT-RHT}
{\bf Y.~F\'elix, S.~Halperin, J.C.~Thomas}, {\it Rational Homotopy
Theory}, Graduate Texts in Mathematics 205, Springer--Verlag
(2001)
  \MR{1802847}

\bibitem{GordonLuecke}
{\bf C. Gordon, J. Luecke},
{\it Knots are determined by their complements},
J. Amer. Math. Soc. 2 (1989) 371--415
  \MR{0965210}


\bibitem{Hovey}
{\bf M.~Hovey}, {\it Model categories}, Mathematical Surveys and
Monographs 63, American Mathematical Society, Providence, RI
(1999)
  \MR{1650134}

\bibitem{Halperin-lecturesminmod}
{\bf S.~Halperin}, {\it Lectures on minimal models}, M\'emoires de
la Soci\'et\'e math\'e\-ma\-tique de France, 230 (1983)
  \MR{0736299}

\bibitem{HalpLemcatDGA}
{\bf S.~Halperin, J.M.~Lemaire}, {\it Notion of category in
differential algebra}, In: Algebraic Topology -- Rational Homotopy,
Springer Lecture Notes in Mathematics 1318,
Berlin--Heidelberg--New--York: Springer--Verlag (1988) 138--154
  \MR{0952577}

\bibitem{Hudson-conc=>iso}
{\bf J.F.P.~Hudson}, {\it Concordance, isotopy, and diffeotopy},
Annals of Math. 91 (1970) 425--448
\MR{0259920}

\bibitem{KahlLVdb}
{\bf T.~Kahl, P. Lambrechts, L.~Vandembroucq}, {\it Mod\`eles de
Quillen des bords homotopiques}, submitted,
available at:\nl
\url{http://www.math.ucl.ac.be/membres/lambrechts/articles.html}

\bibitem{Klein}
{\bf J. Klein}
{\it Poincar\'e duality spaces},
In:  Surveys on surgery theory, Vol. 1, Ann. of Math. Stud. 145, 
Princeton Univ. Press, Princeton, NJ (2000) 135--165
  \MR{1747534}

\bibitem{Klein-2}
{\bf J. Klein}
{\it Poincar\'e duality embeddings and fiberwise homotopy theory, II},
Quart. Jour. Math. Oxford 53 (2002) 319--335 
  \MR{1930266}

\bibitem{Klein-compression}
{\bf J. Klein}
{\it Embedding, compression and fiberwise homotopy theory},
\agtref2{2002}{15}{311}{336}
  \MR{1917055}

\bibitem{Lambrechts-thickening}
{\bf P. Lambrechts}, {\it Cochain models of thickenings and its
application to rational LS--category}, Manuscripta Math. 103
(2000) 143--160
  \MR{1796311}

\bibitem{LS-unstableblowup}
{\bf P. Lambrechts, D.~Stanley}, {\it Examples of rational homotopy types of
blow-ups},
Proc. American Math. Soc. (to appear),
available at:\nl \url{http://www.math.ucl.ac.be/membres/lambrechts/articles.html}

\bibitem{LS-FM2}
{\bf P. Lambrechts, D.~Stanley},
{\it The rational homotopy type of configuration spaces of two points},
Annales Inst. Fourier. 54 (2004) 1029--1052
  \MR{2111020}

\bibitem{LS-stableblowup}
{\bf P. Lambrechts, D.~Stanley},
 {\it The rational homotopy type of a  blow-up in the stable range},
 submitted, available at:\nl
 \url{http://www.math.ucl.ac.be/membres/lambrechts/articles.html}

\bibitem{LSV-TAMS}
{\bf P. Lambrechts, D.~Stanley, L.~Vandembroucq}, {\it Embeddings
up to homotopy type of two-cones in Euclidean space}, Trans.
American Math. Soc. 354 (2002) 3973--4013
  \MR{1926862}

\bibitem{Levitt}
{\bf N. Levitt}, {\it On the structure of Poincar\'e duality spaces},
Topology 7 (1968) 369--388
  \MR{0248839}

\bibitem{Lickorish}
{\bf R. Lickorish}, 
{\it An introduction to knot theory},
Graduate Texts in Mathematics 175, Springer--Verlag (1997)
  \MR{1472978}

\bibitem{Morgan}
{\bf J. Morgan}, {\it The algebraic topology of smooth algebraic
varieties}, Publ. Math. I.H.E.S. {48} (1978) 137--204
  \MR{0516917}

\bibitem{RourkeSanderson}
{\bf C.P.~Rourke, B.J.~Sanderson}, {\it Introduction to
Piecewise-Linear Topology}, Springer Study Edition, Springer--Verlag,  
Berlin--Heidelberg--New-York (1982)
\MR{0665919}

\bibitem{Sullivan}
{\bf D. Sullivan}, {\it Infinitesimal computations in topology},
Publ. Math. I.H.E.S. {47} (1977) 269--331
  \MR{0646078}

\bibitem{Wall-finiteness}
{\bf C.T.C. Wall}, {\it Finiteness conditions for CW--complexes},
Annals of Math.  {81} (1965) 56--69
  \MR{0171284}

\bibitem{Wall-thick}
{\bf C.T.C. Wall}, {\it Classification problems in differential
topology -- IV. Thickenings}, Topology 5 (1966) 73--94
\MR{0192509}

\end{thebibliography}
\end{document}